\documentclass[a4paper]{article}

\usepackage{etex}

\usepackage[english]{babel}
\usepackage[utf8]{inputenc}
\usepackage[T1]{fontenc}
\usepackage{enumerate}
\usepackage{amsmath,amssymb}
\usepackage{mathtools}
\usepackage[amsmath,thmmarks,hyperref]{ntheorem}
\usepackage{dsfont}
\usepackage{color}

\usepackage{tikz,pgfplots}
	\pgfplotsset{compat=newest}
	\pgfplotsset{plot coordinates/math parser=false}
	\newlength\figureheight
	\newlength\figurewidth

\usepackage{geometry}
\usepackage[notref,final]{showkeys}
\usepackage[bookmarksnumbered=true,hidelinks]{hyperref}

\numberwithin{equation}{section}

\theoremnumbering{arabic}
\theoremstyle{break}
\RequirePackage{latexsym}
\theorembodyfont{\itshape}
\theoremsymbol{}
\theoremheaderfont{\normalfont\bfseries}
\theoremseparator{}
	\newtheorem{theorem}{Theorem}[section]
	\newtheorem{lemma}[theorem]{Lemma}
	\newtheorem{corollary}[theorem]{Corollary}

\theorembodyfont{\upshape}
\theoremsymbol{\ensuremath{\clubsuit}}
	\newtheorem{remark}[theorem]{Remark}
	
\theorembodyfont{\upshape}
\theoremsymbol{}

\theoremstyle{nonumberplain}
\theoremheaderfont{\normalfont\bfseries}
\theoremseparator{:}
\theorembodyfont{\normalfont}
\theoremsymbol{\ensuremath{\square}}
\RequirePackage{amssymb}

        \newtheorem{proof}{Proof}
		\qedsymbol{\ensuremath{\square}}

\newcommand{\RR}{\mathds{R}}
\newcommand{\NN}{\mathds{N}}
\newcommand{\eps}{\varepsilon}
\newcommand{\KAK}{\mathcal{K}}
\newcommand{\KAO}{\mathcal{O}}
\newcommand{\KAI}{\mathcal{I}}

\newcommand{\one}{\mathds{1}}
\newcommand{\Bexp}{\mathcal{B}_\mathrm{exp}}
\newcommand{\Bpow}{\mathcal{B}_\mathrm{pow}}
\newcommand{\Minn}{\mathcal{M}_0}
\newcommand{\Mink}{\mathcal{M}_k}
\newcommand{\ua}{{\underline{a}}}
\newcommand{\oa}{{\overline{a}}}
\newcommand{\oI}{{\overline{I}}}
\newcommand{\tnorm}[1]{\left|\!\!\;\left|\!\!\;\left| {#1} \right|\!\!\;\right|\!\!\;\right|}
\newcommand{\norm}[1]{\left\| {#1} \right\|}
\newcommand{\snorm}[1]{\left| {#1} \right|}
\DeclareMathOperator{\sgn}{sgn}
\DeclareMathOperator{\meas}{meas}

\title{Uniform error estimates for general semilinear turning point problems on layer-adapted meshes}
\author{Simon Becher\footnote{Institute of Numerical Mathematics, 
Technical University of Dresden, Dresden D-01062, Germany.
\mbox{e-mail:} Simon.Becher@tu-dresden.de}}
\date{}

\begin{document}

	\maketitle
	
	\begin{abstract}
		We consider a singularly perturbed semilinear boundary value problem of a general
		form that allows various types of turning points. A solution decomposition
		is derived that separates the potential exponential boundary layer terms.
		The problem is discretized using higher order finite elements on suitable 
		constructed layer-adapted meshes. Finally, error estimates uniform with
		respect to the singular perturbation parameter $\varepsilon$ are proven
		in the energy norm.
	\end{abstract}
	
	

\section{Introduction}
Let us consider a singularly perturbed semilinear boundary value problem of the type
\begin{gather}
\label{prob:general}
\begin{aligned}
	- \eps u''(x) + b(x)u'(x) + f(x,u(x)) & = 0, \qquad \text{for } x \in \oI := [\ua,\oa], \\
	u(\ua)= \nu_{-}, \quad u(\oa) & = \nu_+,
\end{aligned}
\end{gather}
where $0 < \eps \ll 1$ and $b,f$ are supposed to be sufficiently smooth.
Furthermore, we assume that there is a continuous function $c$ such that
\begin{gather}
	\label{prob:assNew}
	\partial_u f(x,u) \geq c(x) \geq \gamma > 0, 
	\quad\! \text{for all } (x,u) \in \oI \times \RR,
	\qquad \left( c-\tfrac{1}{2} b'\right)\!(x) \geq \tilde{\gamma} > 0,
	\quad\! \text{for all } x \in \oI.
\end{gather}
A point $\bar{x} \in \oI$ is called turning point of the problem if $b(\bar{x}) = 0$
and for every neighborhood $U$ of $\bar{x}$ there is a point $x \in U \cap \oI$
such that $b(x) \neq 0$. Note that the assumptions in~\eqref{prob:assNew} on $b$, $c$,
and $f$ are very weak and especially allow an arbitrary number, location, and multiplicity
of turning points. But, since these functions are independent of $\eps$, the turning points
are also independent of $\eps$. So, we exclude the situation that an inner turning
point moves to the boundary when $\eps$ goes to zero.

As result of the general setting of problem~\eqref{prob:general} with~\eqref{prob:assNew},
we have to be aware of many (possibly different) layers. One way to treat these layers and
to enable uniform estimates is the use of suitable layer-adapted meshes. This approach
was used by Liseikin in~\cite[Theorem~7.4.2]{Lis01} to prove the uniform first order
convergence of a simple upwind scheme for the considered semilinear problem.

In this paper higher order finite elements shall be analyzed instead. For some
special cases of problem~\eqref{prob:general} it is already known that
optimal order uniform error estimates can be proven on layer-adapted meshes
in an $\eps$-weighted energy norm, see~\cite{RTU15} for linear problems
without turning points or~\cite{Bec16ArXivSunStynes,Bec16ArXivLis} for linear
problems with a single simple attractive interior turning point. But, is
it possible to prove such estimates in our general setting also? And how
could suitable layer-adapted meshes look like?

We shall answer both questions in the following sections. But first more
information about the behavior of the solution, especially about the
appearance and types of layers, is needed. Therefore, a priori estimates
and a solution decomposition are given in Section~\ref{sec:AprioriAndSolDecomp}
together with some comments on the linear version of the considered problem. 
It shows that exponential boundary layers, interior cusp-type layers,
and certain power-type boundary layers could occur. S-type meshes~\cite{RL99}
and the piecewise equidistant meshes proposed by Sun and Stynes~\cite{SS94},
respectively, have proved their worth in handling the first two classes of
layers. Furthermore, it will turn out that the latter grids can be adopted
to the power-type boundary layers by simply adjusting a parameter. So at
the end of Section~\ref{sec:mesh} we are able to give a convenient mesh
construction strategy for the general problem.

The discretization of the problem is presented in Section~\ref{sec:discret}
along with some first notes on the estimation of the error.
Then Section~\ref{sec:errorEst} is devoted to the completing proof of a
uniform error estimate for higher order finite elements.
Several (new) difficulties have to be managed, for example:
\begin{itemize}
	\item The semilinearity of the problem.
	\item The different techniques for the different layers have to be combined.
	\item In general the mesh outside the exponential boundary layer region is not
		quasi uniform and so inverse inequalities, typically used in the analysis of
		S-type meshes, have to be handled with additional care.
	\item Exponential boundary layers of width $\KAO(\sqrt{\eps} \log 1/\sqrt{\eps})$,
		known from reaction-diffusion problems, may also occur when $b \not\equiv 0$.
		Thus, the convection term has to be estimated for such layers also.
	\item The reasoning in the case of a cusp-type layer has to be transferred
		to the case of a power-type boundary layer.
\end{itemize}
Finally, we will have a look at some examples of linear problems with different layers in
Section~\ref{sec:examples}.

Notation: Throughout the paper let $C$ denote a positive generic constant independent
of $\eps$ and the number of mesh intervals $N$. For $S \subset \RR$ we use
the common Sobolev spaces $W^{k,\infty}(S)$, $H^k(S)$, $H_0^1(S)$, and $L^p(S)$.
The spaces of continuously and Lipschitz-continuously differentiable functions
will be written as $C^k(S)$ and $C^{k,1}(S)$, respectively, and used for
$S \subset \RR^2$ also. Furthermore, we shall denote the $L^2$-norm by
$\norm{\cdot}_{0,S}$, the $H^1$-semi norm by $\snorm{\cdot}_{1,S}$, the (essential)
supremum by $\norm{\cdot}_{\infty, S}$, and the $L^p$-norm by $\norm{\cdot}_{L^p(S)}$.
If $S$ is the whole interval, it will be omitted to shorten the notation.

\section{A priori estimates and solution decomposition}
\label{sec:AprioriAndSolDecomp}

A priori estimates for the solution of problem~\eqref{prob:general} with~\eqref{prob:assNew}
can be found, e.g., in~\cite[pp.~73, 74, 95, 96]{Lis01}. We denote by $\Minn := \{\bar{x}_1, \bar{x}_2, \ldots \}$
the set of all points in $(\ua,\oa)$ with $b(\bar{x}_j)=0$ and $b'(\bar{x}_j) < 0$, $j=1,2,\ldots$, i.e.,
all interior turning points, where $b$ changes its sign from $+1$ to $-1$. Note that
$\Minn$ is always finite, see~\cite[p.~73]{Lis01}.

\begin{theorem}[A priori estimates, cf.~\cite{Lis01}]
	\label{th:aPrioriEst}
	Let $q \in \NN$ and suppose that $b \in C^{q}(\oI)$ and $f, f_u \in C^{q}(\oI \times \RR)$.
	Then we have for $k = 0, \ldots, q$ and $x \in \oI$
	\begin{subequations}
	\label{ieq:solDerBound}
	\begin{gather}
	\begin{aligned}
		\big|u^{(k)}(x)\big|
			& \leq C \bigg( 1 + \phi_\ua\big((x-\ua),k,b(\ua),b'(\ua),\eps\big)
				+ \phi_\oa\big((\oa-x),k,-b(\oa),b'(\oa),\eps\big) \\
			& \hspace{18em} + \sum_j \left( \eps^{1/2} + |x-\bar{x}_j| \right)^{\lambda_j-k} \bigg),
	\end{aligned}
	\end{gather}
	where $0 < \lambda_j < c(\bar{x}_j)/|b'(\bar{x}_j)|$ while
	\begin{gather}
		\phi_{\bar{x}}(x,k,a,b,\eps) =
		\begin{cases}
			\eps^{-k} e^{a x/\eps},
				& a < 0, \\
			\eps^{\lambda/2} \left(\eps^{1/2}+x\right)^{-\lambda-k},
				& a = 0, \quad b > 0, \\
			\eps^{-k/2} e^{-\sqrt{(c(\bar{x})+b)} \, x/\sqrt{\eps}},
				& a = 0, \quad 0 \leq -k b < c(\bar{x}) \\
			\left(\eps^{1/2}+x\right)^{\lambda-k} + \eps \left(\eps^{1/2}+x\right)^{-k-2},
				& a = 0, \quad b < 0, \\
			0,	& a > 0,
		\end{cases}
	\end{gather}
	\end{subequations}
	with $0 < \lambda < c(\bar{x})/|b'(\bar{x})|$.
\end{theorem}

Using the standard inequality $1 + x \leq e^x$ we observe that for $\beta > 0$ and $\tilde{x},\tilde{\eps} > 0$
\begin{gather}
	\label{ieq:boundExpAway}
	\tilde{\eps}^{-k} e^{-\beta x/\tilde{\eps}}
		=  \left(\frac{k}{\beta \tilde{x}}\right)^{\!\!k} \left(\frac{\beta \tilde{x}}{k \tilde{\eps}}\right)^{\!\!k}
			e^{-\beta x/\tilde{\eps}}
		\leq \left(\frac{k}{\beta \tilde{x}}\right)^{\!\!k} e^{\beta (\tilde{x}-x)/\tilde{\eps}}
		\leq \left(\frac{k}{\beta \tilde{x}}\right)^{\!\!k}
		\leq C(\tilde{x})
		\qquad \text{when } x \geq \tilde{x}.
\end{gather}
Thus away from the location of boundary and interior layers the solution $u$
and its derivatives can be bounded by a constant. Typically we will use~\eqref{ieq:boundExpAway}
with $\tilde{\eps} \in \{\eps, \sqrt{\eps}\}$.

We want derive a decomposition of the solution $u = S + E$ that separates the
potential exponential boundary layer terms $E$. Unlike for non turning point
problems $S$ is not simply the ``smooth'' part but may consists of power layer terms also.

Let $\Bexp \subseteq \{\ua,\oa\}$ denote the parts of the boundary where
usually an exponential layer occurs. That is
\begin{gather*}
	\Bexp := \big\{ \bar{x} \in \{\ua,\oa\} : b(\bar{x}) \cdot n(\bar{x}) > 0
		\,\text{ or }\, (b(\bar{x}) = b'(\bar{x}) = 0) \big\}
	\quad \text{where} \quad
	n(\bar{x}) := \begin{cases}
		-1,	& \bar{x} = \ua, \\
		1,	& \bar{x} = \oa.
	\end{cases}
\end{gather*}
In an analogous manner let $\Bpow \subseteq \{\ua,\oa\}$, defined by
\begin{gather*}
	\Bpow := \big\{ \bar{x} \in \{\ua,\oa\} : b(\bar{x}) = 0, \, b'(\bar{x}) \neq 0 \big\},
\end{gather*}
denote the parts of the boundary where usually a power-type layer occurs.

For $\bar{x} \in \{\ua,\oa\}$ we define the minimal distance $\delta_{\bar{x}}$
to other (possible) locations of layers by
\begin{gather*}
	\delta_{\bar{x}}
		= \min\!\Big\{ |\bar{x}-y| : y \in \left(\{\ua,\oa\}\setminus \{\bar{x}\}\right)
			\cup \underbrace{\bigcup\nolimits_{j} \{\bar{x}_j\}}_{=\Minn} \Big\}\!.
\end{gather*}
Moreover, for $\bar{x} \in \{\ua,\oa\}$ we define $\bar{\delta}_{\bar{x}}^* \geq 0$
(typical width of a possible exponential layer at $\bar{x}$) by
\begin{gather*}
	\bar{\delta}_{\bar{x}}^* =
	\begin{cases}
		\frac{\eps}{|b(\bar{x})|} \log(1/\eps),
			& b(\bar{x}) \cdot n(\bar{x}) > 0, \\
		\frac{\sqrt{\eps}}{\sqrt{c(\bar{x})}} \log(1/\sqrt{\eps}),
			& b(\bar{x}) = 0, \quad b'(\bar{x}) = 0, \\
		0,	& \text{otherwise}.
	\end{cases}
\end{gather*}

Now, we can prove the following solution decomposition.
\begin{theorem}[Solution decomposition]
	\label{th:solDecomp1}
	Let $q \in \NN$ and suppose that $b \in C^{q}(\oI)$ and $f, f_u \in C^{q}(\oI \times \RR)$.
	Then $u$ has the representation $u = S + E$ with $E= E_\ua + E_\oa$, where for all
	$k = 0, \ldots ,q$ and $x \in \oI$
	\begin{gather}
	\label{ieq:Sbound1}
	\begin{aligned}
		|S^{(k)}(x)|
			& \leq C \bigg( 1 + \phi^S_\ua\big((x-\ua),k,b(\ua),b'(\ua),\eps\big)
				+ \phi^S_\oa\big((\oa-x),k,-b(\oa),b'(\oa),\eps\big) \\
			& \hspace{18em} + \sum_j \left( \eps^{1/2} + |x-\bar{x}_j| \right)^{\lambda_j-k} \bigg),
	\end{aligned}
	\end{gather}
	with $0 < \lambda_j < c(\bar{x}_j)/|b'(\bar{x}_j)|$ and
	\begin{gather}
		\label{ieq:Ebound1}
		|E_\ua^{(k)}(x)|
			\leq C \phi^E_\ua\big((x-\ua),k,b(\ua),b'(\ua),\eps\big),
		\quad
		|E_\oa^{(k)}(x)|
			\leq C \phi^E_\oa\big((\oa-x),k,-b(\oa),b'(\oa),\eps\big).
	\end{gather}
	Here $\phi^S_{\bar{x}}$ and $\phi^E_{\bar{x}}$ are given by
	\begin{align}
		\phi^S_{\bar{x}}(x,k,a,b,\eps) & =
		\begin{cases}
			\eps^{\lambda/2} \left(\eps^{1/2}+x\right)^{-\lambda-k},
				& a = 0, \quad b > 0, \\
			\left(\eps^{1/2}+x\right)^{\lambda-k} + \eps \left(\eps^{1/2}+x\right)^{-k-2},
				& a = 0, \quad b < 0, \\
			0,	& \text{otherwise}
		\end{cases} \label{def:Sphi1}
	\intertext{with $0 < \lambda < c(\bar{x})/|b'(\bar{x})|$ and}
		\phi^E_{\bar{x}}(x,k,a,b,\eps) & =
		\begin{cases}
			\eps^{-k} e^{a x/\eps},
				& a < 0, \\
			\eps^{-k/2} e^{-\sqrt{c(\bar{x})} \, x/\sqrt{\eps}},
				& a = b = 0, \\
			0,	& \text{otherwise}.
		\end{cases} \label{def:Ephi1}
	\end{align}
\end{theorem}

\begin{proof}
	Set $\delta_{\bar{x}}^* = \min\{q \bar{\delta}_{\bar{x}}^*,\delta_{\bar{x}}/2\}$
	for all $\bar{x} \in \{\ua,\oa\}$. For convenience let $\ua^* = \ua + \delta_{\ua}^*$
	and $\oa^* = \oa - \delta_{\oa}^*$. The construction of $\delta_{\bar{x}}^*$ yields
	\begin{gather}
		\label{help_decomp:phiEbound}
		\phi^E_{\bar{x}}\big(\delta_{\bar{x}}^*,0,- b(\bar{x}) \cdot n(\bar{x}),b'(\bar{x}),\eps\big) \geq 
		\phi^E_{\bar{x}}\big(q \bar{\delta}_{\bar{x}}^*,0,- b(\bar{x}) \cdot n(\bar{x}),b'(\bar{x}),\eps\big) =
		\begin{cases}
			\eps^q,		& b(\bar{x}) \cdot n(\bar{x}) > 0, \\
			\eps^{q/2},	& b(\bar{x}) = b'(\bar{x}) = 0, \\
			0,			& \text{otherwise}.
		\end{cases}
	\end{gather}
	This and~\eqref{ieq:boundExpAway}, respectively, gives that up to the $q$th derivative the
	exponential boundary layer terms can be bounded by a constant inside the interval $[\ua^*,\oa^*]$.
	Note that because of $\delta_{\bar{x}}^* \leq \delta_{\bar{x}}/2$ this interval is not empty
	and, moreover, we have
	\begin{gather}
		\label{help_decomp:phiSbound}
		\phi^S_\ua\big((x-\ua),k,b(\ua),b'(\ua),\eps\big)
			+ \phi^S_\oa\big((\oa-x),k,-b(\oa),b'(\oa),\eps\big)
			+ \sum_j \left( \eps^{1/2} + |x-\bar{x}_j| \right)^{\lambda_j-k}
		\leq C
	\end{gather}
	for all $x \in \oI \setminus [\ua^*,\oa^*]$, where the constant $C$
	may depend on $\min\{\delta_{\bar{x}} : \bar{x} \in \Bexp\}$.
	
	In order to prove the decomposition, we adapt an idea from~\cite[p.~23, 24]{RST08}.
	Set $S(x) := u(x)$ for $x \in [\ua^*,\oa^*]$. From~\eqref{ieq:solDerBound},
	\eqref{ieq:boundExpAway}, and~\eqref{help_decomp:phiEbound} we get that $S$
	satisfies~\eqref{ieq:Sbound1} with~\eqref{def:Sphi1} on $[\ua^*,\oa^*]$. Then $S$
	can be extended to a smooth function (i.e., $S \in C^q$) defined on $\oI$ that satisfies
	the requested bound on the whole interval.
	
	We now define $E := u - S$. Obviously, we have $E = 0$ on
	$[\ua^*,\oa^*]$. Combining~\eqref{ieq:solDerBound}, \eqref{ieq:boundExpAway},
	\eqref{help_decomp:phiEbound}, and~\eqref{help_decomp:phiSbound} we get for
	$x \in \oI \setminus [\ua^*,\oa^*]$ that
	\begin{align*}
		\big| E^{(q)}(x) \big|
			& \leq \big| u^{(q)}(x) \big| + \big| S^{(q)}(x) \big| \\
			& \leq C \left(1 + \phi^E_\ua\big((x-\ua),q,b(\ua),b'(\ua),\eps\big)
				+ \phi^E_\oa\big((\oa-x),q,-b(\oa),b'(\oa),\eps\big) \right) \\
			& \leq  C
			\begin{cases}
				\phi^E_\ua\big((x-\ua),q,b(\ua),b'(\ua),\eps\big),	& \text{for } x < \ua^*, \\
				\phi^E_\oa\big((\oa-x),q,-b(\oa),b'(\oa),\eps\big),	& \text{for } \oa^* < x.
			\end{cases}
	\end{align*}
	Integrating $E^{(k)}$ for $k = q, q-1, \ldots, 1$ we inductively gain (recalling that
	$E = 0$ on $[\ua^*,\oa^*]$)
	\begin{align*}
		\big| E^{(k-1)}(x) \big|
			\leq \left| \int_{\ua^*}^x E^{(k)}(s) \, ds \right|
		  & \leq C \int_{\ua^*}^x \phi^E_\ua\big((x-\ua),k,b(\ua),b'(\ua),\eps\big) \, ds \\
		  & \leq C \phi^E_\ua\big((x-\ua),k-1,b(\ua),b'(\ua),\eps\big),
		  	\qquad \text{for } x < \ua^*
	\end{align*}
	and analogously
	\begin{gather*}
		\big| E^{(k-1)}(x) \big| \leq C \phi^E_\oa\big((\oa-x),k-1,-b(\oa),b'(\oa),\eps\big),	
			\qquad \text{for } \oa^* < x.
	\end{gather*}
	Since $E = 0$ on $[\ua^*,\oa^*] \neq \emptyset$ by construction, we simply have to set
	\begin{gather*}
		E_\ua =
			\begin{cases}
				E, & \text{on } [\ua,\ua^*], \\
				0, & \text{otherwise},
			\end{cases}
		\qquad \text{and} \qquad
		E_\oa =
			\begin{cases}
				E, & \text{on } [\oa^*,\oa], \\
				0, & \text{otherwise},
			\end{cases}
	\end{gather*}
	respectively, to get the two terms $E_\ua$ and $E_\oa$. Thus the statement is proven.
\end{proof}

\begin{remark}[Comments on the linear version of the problem]
	Consider the linear version of problem~\eqref{prob:general} where $f(x,u(x))$ is
	replaced by $c(x) u(x) - f(x)$, i.e.,
	\begin{gather}
	\label{prob:generalLin}
	\begin{aligned}
		Lu := - \eps u''(x) + b(x)u'(x) + c(x)u(x) & = f(x), \qquad \text{for } x \in \oI := [\ua,\oa], \\
		u(\ua)= \nu_{-}, \quad u(\oa) & = \nu_+,
	\end{aligned}
	\end{gather}
	with $0 < \eps \ll 1$, problem data $b,c,f$ sufficiently smooth, and suppose that
	\begin{subequations}
	\label{prob:assOld}
	\begin{align}
		\text{for all } x \in \oI: \qquad b(x) = 0 \quad & \Longrightarrow \quad
			c(x)>0, \label{prob:ass1} \\
		\text{for all } x \in \oI: \qquad b(x) = 0 \quad & \Longrightarrow \quad
			\!\left(c - \tfrac{1}{2}b'\right)(x) >0. \label{prob:ass2}
	\end{align}
	\end{subequations}
	Then we may assume without loss of generality (for $\eps$ sufficiently small)
	that~\eqref{prob:assNew} holds, i.e.,
	\begin{gather*}
		c(x) \geq \gamma > 0, \qquad \left( c-\tfrac{1}{2} b'\right)\!(x) \geq \tilde{\gamma} > 0,
		\qquad \text{for all } x \in \oI.
	\end{gather*}
	This can always be achieved by a suitable problem transformation,
	see Appendix~\ref{app:probTrans} for a proof of this statement. Such transformations
	are widely known in the case that $b(x) \neq 0$ and was also studied in~\cite{SS94} in
	the case of a single interior turning point. But to the authors knowledge this
	statement is new in this general setting.
	
	Moreover, we want to note that for the linear problem a solution decomposition
	similar to that of Theorem~\ref{th:solDecomp1} can also be derived using boundary
	layer corrections. The interested reader is referred to Appendix~\ref{app:solDecomp}
	for more details.
\end{remark}

\section{FEM discretization}
\label{sec:discret}

In order to fix the notation we want to present the discretization of the problem
by higher order finite elements now. We will consider homogeneous Dirichlet
boundary conditions $\nu_{-} = \nu_{+} = 0$ only. Note that these can be easily ensured
by modifying the nonlinear term $f(x,u)$. Indeed, set
$\tilde{u}(x) = \left((\oa-x) \nu_{-} + (x-\ua) \nu_{+}\right)/(\oa-\ua)$.
Then $u + \tilde{u}$ solves~\eqref{prob:general} when $u$ itself solves
\begin{gather*}
	- \eps u''(x) + b(x)u'(x) + \tilde{f}(x,u(x)) = 0, \qquad \text{for } x \in \oI := [\ua,\oa],
	\qquad \quad
	u(\ua)= u(\oa) = 0,
\end{gather*}
where $\tilde{f}(x,u) = f(x,u+\tilde{u}(x)) - \eps \tilde{u}''(x) + b(x) \tilde{u}'(x)$
for $u \in \RR$. Obviously, the assumption~\eqref{prob:assNew} also holds for $\tilde{f}$.
	
For $v, w \in V := H_0^1((\ua,\oa))$ we set
\begin{gather*}
	B_\eps\!\left(v,w\right) := \left(\eps v',w'\right) + \left(b v', w\right) + \left(f(\cdot,v),w\right)\!.
\end{gather*}
Then we obtain the following weak formulation of~\eqref{prob:general} with $\nu_{-} = \nu_{+} = 0$:
\medskip

Find $u \in V$ such that
\begin{gather}
	\label{weakprob:general}
	B_\eps\!\left(u,v\right) = 0, \qquad \text{ for all } v \in V.
\end{gather}

By $\ua = x_{0} < \ldots < x_i < \ldots < x_N = \oa$ an arbitrary mesh is given
on the interval $[\ua,\oa]$. Let $h_i := x_i-x_{i-1}$ denote the mesh interval
lengths. We define the trial and test space $V^N$ by
\begin{gather*}
	V^N := \left\{ v \in C([\ua,\oa]) : v|_{(x_{i-1},x_i)} \in P_k((x_{i-1},x_i))\, \forall i, \, v(-1) = v(1) = 0 \right\}
		\subset V
\end{gather*}
where $k \geq 1$. The space $P_k((x_a,x_b))$ comprises all polynomials up to order $k$
over $(x_a,x_b)$. The discrete problem arises from replacing $V$
in~\eqref{weakprob:general} by the finite dimensional subspace $V^N$: \medskip

Find $u^N \in V^N$ such that
\begin{gather}
	\label{dprob:general}
	B_\eps\!\left(u^N,v^N\right) = 0, \qquad \text{ for all } v^N \in V^N.
\end{gather}

Let $v^I$ denote the standard Lagrangian interpolant into $V^N$ of $v \in V$.
As interpolation points we choose the mesh points and $k-1$ (arbitrary) inner points per interval.
For example uniform or Gau{\ss}-Lobatto points could be used.

Assuming $v \in W^{k+1,\infty}((x_{i-1},x_i))$,
the standard interpolation theory leads to the error estimates: 
For all $j = 0, \ldots, k+1$
\begin{gather}
	\norm{(v-v^I)^{(j)}}_{\infty,(x_{i-1},x_i)}
		\leq C h_i^{k+1-j} \norm{v^{(k+1)}}_{\infty,(x_{i-1},x_i)} \label{standInt} \\
	\intertext{and}
	\norm{(v-v^I)^{(j)}}_{\infty,(x_{i-1},x_i)}
		\leq C \bigl\| v^{(j)} \bigr\|_{\infty,(x_{i-1},x_i)} \label{standIntInfty}.
\end{gather}
The constant $C$ depends on the location of the inner interpolation points.
Furthermore, for all $v^N \in V^N$ the inverse inequality
\begin{gather}
	\label{invIneq}
	\norm{(v^N)'}_{L^p((x_{i-1},x_i))} \leq C h_i^{-1} \norm{v^N}_{L^p((x_{i-1},x_i))}
\end{gather}
holds for $p \in [1,\infty]$.

Now we have a closer look on the nonlinear term. The following identity
will be exploited several times later on. For a function $g = g(x,v)$ with
$g, g_v \in C(\oI \times \RR)$ we have
\begin{gather}
	\label{eq:nonlinIdentity}
	g(x,v_1)-g(x,v_2) = \left[\int_0^1 \partial_v g(x,v_2+s(v_1-v_2)) \,ds\right] (v_1-v_2).
\end{gather}
Especially, choosing $g(x,v) = f(x,v)-\frac{1}{2}b'v$ and noting that
\begin{gather*}
	\partial_v g(x,v)
		= \partial_u f(x,v) - \tfrac{1}{2}b'(x)
		\geq c(x) - \tfrac{1}{2}b'(x)
		\geq \tilde{\gamma} > 0
\end{gather*}
by~\eqref{prob:assNew}, we obtain for $v_1, v_2 \in V$ and shortly $v = v_1 - v_2$
\begin{gather}
\label{ieq:semiCoercive}
\begin{aligned}
	B_\eps\!\left(v_1,v\right)-B_\eps\!\left(v_2,v\right)
		& = \eps \left(v',v'\right) + \left(b v', v\right) + \left(f(\cdot,v_1)-f(\cdot,v_2),v\right) \\
		& = \eps \snorm{v}_1^2 + \left([f(\cdot,v_1) - \tfrac{1}{2}b' v_1] - [f(\cdot,v_2)-\tfrac{1}{2}b' v_2],v\right) \\
		& \geq \eps \snorm{v}_1^2 + \left((c-\tfrac{1}{2}b')v,v\right) \\
		& \geq \eps \snorm{v}_1^2 + \tilde{\gamma} \norm{v}_0^2.
\end{aligned}
\end{gather}
Therefore, for the analysis of the problem an appropriated weighted energy norm
$\tnorm{\cdot}_\eps$ is given by
\begin{gather*}
	\tnorm{v}_\eps := \left( \eps \snorm{v}_1^2 + \tilde{\gamma} \norm{v}_0^2\right)^{1/2}.
\end{gather*}
Note that~\eqref{ieq:semiCoercive} implies the uniqueness of the weak and the discrete solution.
In the case of a linear problem $B_\eps\!\left(\cdot,\cdot\right)$ is a bilinear
form which is uniformly coercive over $H_0^1((\ua,\oa)) \times H_0^1((\ua,\oa))$
with respect to $\tnorm{\cdot}_\eps$ due to~\eqref{ieq:semiCoercive}.

As usual, to estimate the discretization error a splitting is used, that is
\begin{gather}
	\label{splitError}
	u-u^N = (u-u^I) + (u^I-u^N).
\end{gather}
Combining~\eqref{ieq:semiCoercive} and the problem formulations~\eqref{weakprob:general}
and~\eqref{dprob:general}, we conclude
\begin{gather}
	\label{startEstimation}
	\tnorm{u^N-u^I}_{\eps}^2
		\leq B_{\eps}(u^N,v^N) - B_{\eps}(u^I,v^N)
		= B_{\eps}(u,v^N) - B_{\eps}(u^I,v^N)
\end{gather}
where $v^N = u^N-u^I \in V^N \subset V$. Unlike for linear problems~\cite{Bec16ArXivSunStynes,RTU15}
the right hand side of~\eqref{startEstimation} contains nonlinear terms in $u$ and $u^I$.
So the nonlinearity needs additional consideration.

Fortunately, from inverse monotonicity properties,
see~\cite[p.~74]{Lis01}, we have for the solution $u$ of problem~\eqref{prob:general} that
\begin{gather*}
	\norm{u}_\infty \leq \max \big\{|\nu_{-}|,|\nu_{+}|, \norm{f(\cdot,0)/c(\cdot)}_\infty \big\}.
\end{gather*}
Hence, with~\eqref{standIntInfty} we get
\begin{gather*}
	\norm{u}_\infty + \norm{u^I}_\infty \leq \tilde{C},
\end{gather*}
where $\tilde{C}>0$ is independent of $\eps$ and depends on the problem data and
the location of the inner interpolation points only. Furthermore, for
$f_u \in C(\oI \times \RR)$ we can define
\begin{gather*}
	\tilde{C}_f := \max\big\{\left|\partial_u f(x,v)\right| : (x,v) \in \oI \times [-\tilde{C},\tilde{C}] \big\}.
\end{gather*}
Therefore, using~\eqref{eq:nonlinIdentity} with $g(x,v) = f(x,v)$ yields
\begin{gather}
	\label{ieq:nonlinToLin}
	\left| f\big(x,u(x)\big)-f\big(x,u^I(x)\big) \right|
		\leq \underbrace{\left[\int_0^1 \left|\partial_u f\big(x,u^I(x)+s\big(u(x)-u^I(x)\big)\big) \right|
				ds\right]}_{\leq \tilde{C}_f}
			\left|u(x)-u^I(x) \right|
\end{gather}
which simply reduces the difference of the nonlinear terms to a basic interpolation error. Thus,
the nonlinearity raises no further difficulties in the error estimation.

Under the above assumptions the existence of unique solutions for~\eqref{weakprob:general}
and~\eqref{dprob:general} is clear. For details we refer to~\cite[Section~3.5]{GRS07}.

\section{Layer-adapted meshes}
\label{sec:mesh}

In this section we want to describe how a layer-adapted mesh
for problem~\eqref{prob:general} with~\eqref{prob:assNew} could
be constructed, see Section~\ref{subsec:genMesh}. Before some
layer-adapted meshes are presented each of which can be used to
capture a single layer-type. In order to describe these meshes we
shall always consider the interval $[0,1]$ and assume that the
layer is located at zero.

\subsection{S-type meshes for exponential layers}
\label{subsec:S-typeMesh}

S-type meshes are very popular in the cases of exponential layers. They are
characterized by a very fine mesh in the layer region and a coarse mesh
away from the layer. Dependent on the detailed structure of the exponential
layer term $\phi_E$ the transition point $\tau$ from the fine to the coarse
part is defined by
\begin{gather}
	\label{def:tau}
	\tau = \rho \frac{\tilde{\eps}}{\beta} \ln N, \qquad
	\tilde{\eps} = 
	\begin{cases}
		\eps,
			& \text{if } \big|\phi_E^{(j)}(x)\big| \leq C \eps^{-j} e^{-\beta x/\eps}, \quad 0 \leq j \leq k+1, \\
		\sqrt{\eps},
			& \text{if } \big|\phi_E^{(j)}(x)\big| \leq C \eps^{-j/2} e^{-\beta x/\sqrt{\eps}}, \quad 0 \leq j \leq k+1,
	\end{cases}
\end{gather}
with $\rho > 0$. This definition yields $\big| \phi_E(\tau) \big| \leq C N^{-\rho}$.

Now, we set $x_{N/2} = \tau$ and define the mesh on the interval $[0, x_{N/2}]$ by
a mesh generating function $\varphi$, where we suppose that $\varphi$ is continuous,
monotonically increasing, piecewise continuously differentiable, and furthermore
satisfies $\varphi(0)=0$ and $\varphi(\frac{1}{2}) = \ln N$. On the interval
$[x_{N/2},1]$ a coarse mesh is generated. Often an equidistant
partition (with $h_i \leq C N^{-1}$) is chosen. In this case we have
\begin{gather*}
	x_i = 
	\begin{cases}
		\rho \frac{\tilde{\eps}}{\beta} \varphi(t_i),
			& t_i =i/N, \, i = 0, 1,\ldots, N/2, \\
		1- (1-\tau) \frac{2(N-i)}{N},
			& i = N/2+1, \ldots, N.
	\end{cases}
\end{gather*}
Moreover, the mesh characterizing function $\psi$ is defined by
\begin{gather*}
	\varphi = - \ln \psi \qquad
	\left( \Rightarrow \; \psi(t_i) = \exp\left(-\frac{\beta x_i}{\rho\tilde{\eps}}\right)\right)\!.
\end{gather*}

The further analysis is based on the assumption (typical for FEMs on S-type meshes)
\begin{gather}
	\label{ass:meshGenFuncDerBound}
	\max \varphi' \leq C N.
\end{gather}
Therefore, by the construction of the mesh we get for $1 \leq i \leq N/2$
\begin{gather}
	\label{ieq:hiBoundTildeEps}
	h_i = x_i - x_{i-1}
		= \rho \frac{\tilde{\eps}}{\beta} \big(\varphi(t_i) - \varphi(t_{i-1})\big)
		\leq \rho \frac{\tilde{\eps}}{\beta} N^{-1} \max\varphi'
		\leq C \tilde{\eps}.
\end{gather}

We now prove some estimates for the interpolation error.
\begin{lemma}
	\label{le:IntErrorS-type}
	Under the assumption~\eqref{ass:meshGenFuncDerBound} and for $\rho \geq k+1$ we have
	\begin{align*}
		\norm{\phi_E - \phi_E^I}_{\infty,[0,\tau]}
			& \leq C \left(N^{-1} \max |\psi'|\right)^{k+1},
		& \norm{\phi_E - \phi_E^I}_{\infty,[\tau,1]}
			& \leq C N^{-\rho},\\		
		\tilde{\eps}^{1/2} \snorm{\phi_E - \phi_E^I}_{1,[0,\tau]}
			& \leq C \left(N^{-1} \max |\psi'|\right)^k,
		& \tilde{\eps}^{1/2} \snorm{\phi_E - \phi_E^I}_{1,[\tau,1]}
			& \leq C N^{-\rho}, \\
		\norm{\phi_E - \phi_E^I}_{0,[0,\tau]}
			& \leq C \tilde{\eps}^{1/2} (\ln N)^{1/2} \left(N^{-1} \max |\psi'|\right)^{k+1},
		& \norm{\phi_E - \phi_E^I}_{0,[\tau,1]}
			& \leq C N^{-\rho},
	\intertext{or if $\rho \geq k+3/2$} 
		\norm{\phi_E - \phi_E^I}_{0,[0,\tau]}
			& \leq C \tilde{\eps}^{1/2} \left(N^{-1} \max |\psi'|\right)^{k+1}.
	\end{align*}
	Here $\phi_E^I \in V^N$ denotes a standard Lagrangian interpolant of $\phi_E$
	($k+1$ points per mesh interval).
\end{lemma}
\begin{proof}
	The proof follows~\cite[Section~2]{RTU15} but is generalized with respect
	to the different width of the exponential layer. We first study the
	interval $[0,\tau]$, so let $1 \leq i \leq N/2$. Applying some standard
	interpolation error estimates we obtain
	\begin{gather*}
		\norm{\phi_E-\phi_E^I}_{\infty,[x_{i-1},x_i]}
			\leq C h_i^{k+1} \norm{\phi_E^{(k+1)}}_{\infty,[x_{i-1},x_i]}\!.
	\end{gather*}
	Together with
	\begin{gather}
		\label{ieq:hiFine}
		h_i \leq \rho \frac{\tilde{\eps}}{\beta} N^{-1} \max \varphi'
			\leq \rho \frac{\tilde{\eps}}{\beta} N^{-1} \max | \psi'| e^{\beta x_i/(\rho \tilde{\eps})}
	\end{gather}
	and the general bound $h_i \leq C \tilde{\eps}$ (from~\eqref{ieq:hiBoundTildeEps})
	this yields for $\rho \geq k+1$
	\begin{gather*}
		\norm{\phi_E-\phi_E^I}_{\infty,[x_{i-1},x_i]}
			\leq C \left(N^{-1} \max |\psi'|\right)^{k+1} e^{\beta \left(\frac{(k+1)}{\rho}x_i-x_{i-1}\right)/\tilde{\eps}}
			\leq C \left(N^{-1} \max |\psi'|\right)^{k+1}.
	\end{gather*}
	
	The first bound in the $L^2$-norm follows from the maximum norm estimate since
	\begin{gather*}
		\norm{\phi_E-\phi_E^I}_{0,[0,\tau]}
			\leq \tau^{1/2} \norm{\phi_E-\phi_E^I}_{\infty,[0,\tau]}
			\leq C \tilde{\eps}^{1/2} (\ln N)^{1/2} \left(N^{-1} \max |\psi'|\right)^{k+1}.
	\end{gather*}
	
	In the integral-based norms (with $j = 0,1$) we can also proceed as follows
	\begin{align*}
		\snorm{\phi_E-\phi_E^I}_{1-j,[x_{i-1},x_i]}^2
		  & \leq C h_i^{2(k+j)} \int_{x_{i-1}}^{x_i} \!\left(\phi_E^{(k+1)} \right)^2
			\leq C h_i^{2(k+j)} \tilde{\eps}^{-(2k+1)}
				\left(e^{-2\beta x_{i-1}/\tilde{\eps}}-e^{-2\beta x_i/\tilde{\eps}}\right) \\
		  & \leq C \tilde{\eps}^{-1+2j} \left(N^{-1} \max |\psi'|\right)^{2(k+j)}
		  		e^{2(k+j) \beta x_i/(\rho\tilde{\eps})} e^{-2\beta x_{i-1/2}/\tilde{\eps}} \sinh \frac{\beta h_i}{\tilde{\eps}}.
	\end{align*}
	By~\eqref{ieq:hiBoundTildeEps} we have $\sinh \frac{\beta h_i}{\tilde{\eps}} \leq \sinh C \beta \leq C$ and thus
	\begin{gather*}
		\sinh \frac{\beta h_i}{\tilde{\eps}}
			\leq C \frac{\beta h_i}{\tilde{\eps}}
			= C \int_{t_{i-1}}^{t_i} \varphi'
			= C \int_{t_{i-1}}^{t_i} \left(\frac{-\psi'}{\psi}\right)
			\leq C e^{\beta x_i/(\rho \tilde{\eps})} \int_{t_{i-1}}^{t_i} (-\psi').
	\end{gather*}
	Combining
	\begin{gather*}
		e^{2(k+j) \beta x_i/(\rho\tilde{\eps})} e^{-2\beta x_{i-1/2}/\tilde{\eps}}e^{\beta x_i/(\rho \tilde{\eps})}
			= e^{2\beta\left(\frac{(k+j+1/2)}{\rho}x_i-x_{i-1/2}\right)/\tilde{\eps}}
			\leq C
	\end{gather*}
	which holds for $\rho \geq k+j+1/2$ and
	\begin{gather*}
		\int_0^{1/2} (-\psi') = \psi(0) - \psi(1/2) \leq 1
	\end{gather*}
	we get
	\begin{gather*}
		\tilde{\eps}^{1-2j} \snorm{\phi_E-\phi_E^I}_{1-j,[0,\tau]}^2
			\leq C \left(N^{-1} \max |\psi'|\right)^{2(k+j)}.
	\end{gather*}
	
	The estimate in the interval $[\tau, 1]$ is based on the smallness of the
	exponential boundary layer term. Using some stability properties of the
	interpolant, we conclude from~\eqref{def:tau}
	\begin{gather*}
		\norm{\phi_E-\phi_E^I}_{\infty,[\tau,1]}
			\leq \norm{\phi_E}_{\infty,[\tau,1]} + \norm{\phi_E^I}_{\infty,[\tau,1]}
			\leq C \norm{\phi_E}_{\infty,[\tau,1]}
			\leq C \snorm{\phi_E(\tau)}
			\leq C N^{-\rho}.
	\end{gather*}
	Immediately, this implies
	\begin{gather*}
		\norm{\phi_E-\phi_E^I}_{0,[\tau,1]}
			\leq C \norm{\phi_E-\phi_E^I}_{\infty,[\tau,1]}
			\leq C N^{-\rho}.
	\end{gather*}
	Using the stability/error estimate for the interpolator we also get
	\begin{align*}
		\snorm{\phi_E-\phi_E^I}_{1,[\tau,1]}^2
			\leq C \snorm{\phi_E}_{1,[\tau,1]}^2
			\leq C \int_\tau^1 \tilde{\eps}^{-2} e^{-2\beta x/\tilde{\eps}}
			\leq C \tilde{\eps}^{-1} \left(e^{-2\beta\tau/\tilde{\eps}} - e^{-2\beta/\tilde{\eps}}\right)
			\leq C \tilde{\eps}^{-1} N^{-2\rho}
	\end{align*}
	which completes the proof.
\end{proof}

\begin{remark}
	\label{rem:coarseMeshPart}
	The estimates of Lemma~\ref{le:IntErrorS-type} on the coarse part $[\tau,1]$
	are independent of the exact structure of the mesh in this interval. The proof
	only uses the smallness of $\phi_E$, i.e., $\big| \phi_E(x) \big| \leq C N^{-\rho}$
	for $x \geq \tau$.
\end{remark}

\begin{lemma}
	\label{le:xTimesIntErrorS-type}
	Under the assumption~\eqref{ass:meshGenFuncDerBound} and for $\rho \geq k+1$, $\ell \geq 1$ we have
	\begin{gather*}
		\norm{x^\ell(\phi_E - \phi_E^I)'}_{\infty,[0,\tau]}
			\leq C (2\tilde{\eps})^{\ell-1} \left(N^{-1}\max|\psi'|\right)^{\!k}.
	\end{gather*}
\end{lemma}
\begin{proof}
	Let $1 \leq i \leq N/2$. Combining standard interpolation estimates,
	the bounds for $\phi_E$, and~\eqref{ieq:hiFine} we obtain
	\begin{align*}
		&\norm{x^\ell(\phi_E-\phi_E^I)'}_{\infty,[x_{i-1},x_i]} \\
		& \quad \leq C (x_{i-1}+h_i)^\ell h_i^k \norm{\phi_E^{(k+1)}}_{\infty,[x_{i-1},x_i]}
			\leq C 2^{\ell-1} \left(x_{i-1}^\ell + h_i^\ell\right)\left(\frac{h_i}{\tilde{\eps}}\right)^{\!\!k}
				\tilde{\eps}^{-1} e^{-\beta x_{i-1}/\tilde{\eps}} \\
		& \quad \leq C (2\tilde{\eps})^{\ell-1} \left(\left(N^{-1}\max|\psi'|\right)^{\!k}
				\left(\frac{x_{i-1}}{\tilde{\eps}}\right)^{\!\ell} e^{\beta k x_i/(\rho\tilde{\eps})}
			+ \left(N^{-1}\max|\psi'|\right)^{\!k+1} e^{\beta (k+1) x_i/(\rho\tilde{\eps})} \right) 
				e^{-\beta x_{i-1}/\tilde{\eps}}.
	\end{align*}
	From the well known inequality $1 + x \leq e^x$ we conclude
	\begin{gather*}
		\left(\frac{x_{i-1}}{\tilde{\eps}}\right)^{\!\ell} e^{- \beta x_{i-1}/\tilde{\eps}}
			\leq \left(\frac{\ell(k+1)}{\beta}\right)^{\!\ell}
				e^{- \beta k x_{i-1}/((k+1)\tilde{\eps})}.
	\end{gather*}
	Hence, we have for $\rho \geq k+1$
	\begin{align*}
		&\norm{x(\phi_E-\phi_E^I)'}_{\infty,[x_{i-1},x_i]} \\
		& \qquad \leq C (2\tilde{\eps})^{\ell-1} \left(\left(N^{-1}\max|\psi'|\right)^{\!k}
				e^{\beta \left(\frac{k}{\rho} x_i -\frac{k}{k+1} x_{i-1}\right)/\tilde{\eps}}
			+ \left(N^{-1}\max|\psi'|\right)^{\!k+1}
				e^{\beta \left(\frac{k+1}{\rho} x_i-x_{i-1}\right)/\tilde{\eps}} \right) \\
		& \qquad \leq C (2\tilde{\eps})^{\ell-1} \left(\left(N^{-1}\max|\psi'|\right)^{\!k}
				+ \left(N^{-1}\max|\psi'|\right)^{\!k+1} \right)\!,
	\end{align*}
	where $x_i - x_{i-1} = h_i \leq C \tilde{\eps}$ was used in the last inequality.
\end{proof}

\subsection{Piecewise equidistant mesh for power-type layers}
\label{subsec:SunStynesMesh}

The piecewise equidistant mesh presented in this section was first introduced
by Sun and Stynes~\cite[Section~5.1]{SS94} to treat an interior cusp-type layer.
We will see that it can be easily adapted such that uniform estimates are possible
for all kinds of power-type layers appearing in~\eqref{def:Sphi1} as well. This
fact is based on the observation that all of these layer terms can be bounded as
\begin{gather}
	\label{ieq:phiS}
	\snorm{\phi_S^{(j)}(x)} \leq C \left( 1 + \left(\eps^{1/2} + x \right)^{\!\lambda-j}\right)\!,
	\quad 0 \leq j \leq k+1,
\end{gather}
with $\lambda \geq 0$.

We shall assume in the following that $\lambda \in [0,k+1)$. This is the most
difficult case since otherwise all crucial derivatives of $\phi_S$ could be
bounded by a constant independent of $\eps$ which would allow uniform estimates
for the layer terms using standard methods on uniform meshes.

The mesh parameters are determined as in~\cite[Section~3]{Bec16ArXivSunStynes}.
For $\eps \in (0,1]$ and given positive integer $N$ we set
\begin{gather}
	\label{sigma}
	\sigma = \max\big\{\eps^{(1-\lambda/(k+1))/2}, N^{-(2k+1)} \big\}
\end{gather}
and
\begin{gather}
	\label{KAK}
	\KAK = \left\lfloor 1 - \frac{\ln(\sigma)}{\ln(10)} \right\rfloor,
\end{gather}
where $\lfloor z \rfloor$ denotes the largest integer less or equal to $z$.

The piecewise equidistant mesh is constructed in two steps: First, we divide the
interval $(0,1]$ into the $\KAK+1$ subintervals
$(0,10^{-\KAK}], \, (10^{-\KAK},10^{-\KAK+1}], \ldots , (10^{-1},1]$. Afterwards
each of these subintervals is partitioned uniformly into $\lfloor N/(\KAK+1) \rfloor$
parts, where for simplicity we assume that $\lfloor N/(\KAK+1) \rfloor = N/(\KAK+1)$.
Thus, by construction we have
\begin{gather}
	\label{eq:hiKAK}
	h_i = (\KAK+1)10^{-\KAK}N^{-1}, \quad \text{for} \quad x_i \in (0,10^{-\KAK}]
\end{gather}
and
\begin{gather}
	\label{eq:hiEll}
	h_i = 9(\KAK+1)10^{-l}N^{-1}, \quad \text{for} \quad x_i \in (10^{-l},10^{-l+1}]
		\quad \text{and} \quad l = 1,\ldots,\KAK.
\end{gather}

The properties of the logarithm together with~\eqref{sigma} and~\eqref{KAK} yield
\begin{gather*}
	\KAK + 1 \leq 2 - \frac{\ln(\sigma)}{\ln(10)}
		\leq 2 + \min\left\{\frac{1-\lambda/(k+1)}{2}\frac{|\ln(\eps)|}{\ln(10)},(2k+1)\frac{\ln(N)}{\ln(10)}\right\}.
\end{gather*}
This ensures for sufficiently large $N$ (dependent on $k$) that the number of subintervals
$\KAK + 1$ is less than the number of mesh intervals $N$ since
\begin{gather}
	\label{ieq:KAKtoLN}
	\KAK + 1 \leq C \ln N.
\end{gather}
Moreover, from~\eqref{KAK} we easily see that
\begin{gather}
	\label{ieq:sigma}
	10^{-1} \sigma \leq 10^{-\KAK} < \sigma.
\end{gather}
Note that by a simple modification of the construction one can guarantee
that the mesh consists of exactly $N$ mesh intervals, see~\cite[Section~6]{SS94}.
Let $\Phi_{N,\lambda}$ denote the associated mesh generating function which
is continuous and piecewise linear.

The next lemma is taken from~\cite{Bec16ArXivSunStynes}. Note that the arguments
used therein also works for the (formal) choice $\lambda = 0$.
\begin{lemma}[{see~\cite[Lemma~3.1]{Bec16ArXivSunStynes}}]
	\label{le:sunStynesh^kBounds}
	Let $j=0,1$. The following inequalities hold
	\begin{align}
		h_i^{k+1-j}\left(x_{i-1}+\eps^{1/2}\right)^{\lambda-(k+1-j)}
			& \leq C \left((\KAK + 1)N^{-1}\right)^{k+1-j}, \quad \text{for} \quad x_i \in (10^{-\KAK},1], \\
		h_i^{k+1-j}\left(x_{i-1}+\eps^{1/2}\right)^{\lambda-(k+1-j)}
			& \leq C\left(i-1\right)^{-(k+1-j)}, \quad \text{for} \quad x_i \in (x_1,10^{-\KAK}]. \label{ieq:sunStynes_im1}
	\intertext{If $\sigma = \eps^{(1-\lambda/(k+1))/2}$, then}
		h_i^{k+1-j}\left(x_{i-1}+\eps^{1/2}\right)^{\lambda-(k+1-j)}
			& \leq C \left((\KAK + 1)N^{-1}\right)^{k+1-j}, \quad \text{for} \quad x_i \in (0,10^{-\KAK}].
	\end{align}
	In general, the mesh interval length can be bounded by
	\begin{gather*}
		h_i \leq (\KAK + 1)N^{-1}.
	\end{gather*}
	Furthermore, in the case of $\sigma = N^{-(2k+1)}$, we have
	\begin{gather}
		\label{ieq:sunStynesX_1}
		x_1 = h_1 \leq (\KAK+1) N^{-2(k+1)}.
	\end{gather}
\end{lemma}

Using the techniques of~\cite[Lemma~3.2]{Bec16ArXivSunStynes} (also here the reasoning
can be adopted for the (formal) choice $\lambda = 0$), we obtain the following
interpolation error estimates on the layer-adapted piecewise equidistant mesh.
\begin{lemma}[{cf.~\cite[Lemma~3.2]{Bec16ArXivSunStynes}}]
	\label{le:IntErrorSunStynes}
	Let $\phi_S$ satisfy~\eqref{ieq:phiS} and let $\phi_S^I \in V^N$ be its interpolant
	on the piecewise equidistant mesh given by~\eqref{sigma} -- \eqref{eq:hiEll}. Then
	\begin{gather}
		\label{ieq:L2IntSunStynes}
		\norm{\phi_S-\phi_S^I}_0 \leq C \left((\KAK + 1) N^{-1}\right)^{k+1}
	\end{gather}
	and
	\begin{gather}
		\label{ieq:energyIntSunStynes}
		\tnorm{\phi_S-\phi_S^I}_\eps + \norm{x(\phi_S-\phi_S^I)'}_0 \leq C \left((\KAK + 1) N^{-1}\right)^{k}.
	\end{gather}
\end{lemma}

\begin{remark}
	\label{rem:xihiInv}
	By the construction of the mesh we have
	\begin{gather*}
		x_i h_i^{-1} \leq C N.
	\end{gather*}
	Indeed, for $x_i \in (0,10^{-\KAK}]$ it holds
	\begin{gather*}
		x_i h_i^{-1} \leq \frac{10^{-\KAK}}{(\KAK + 1) 10^{-\KAK} N^{-1}} \leq N (\KAK+1)^{-1}
	\end{gather*}
	and for $x_i \in (10^{-l},10^{-l+1}]$ with $l = 1,\ldots,\KAK$
	\begin{gather*}
		x_i h_i^{-1} \leq \frac{10^{-l+1}}{9(\KAK + 1) 10^{-l} N^{-1}} \leq \frac{10}{9} N (\KAK+1)^{-1}.
	\end{gather*}
\end{remark}

\subsection{Layer-adapted mesh for general turning point problem}
\label{subsec:genMesh}

First of all we recall some relevant observations of the previous sections:
\begin{itemize}
	\item For power-type layers of the form~\eqref{ieq:phiS} an adaption of the mesh
		is only needed when $\lambda < k+1$. Such a situation is present in the vicinity
		of $\bar{x} \in \Mink \cup \Bpow$ where
		\begin{gather*} 
			\Mink := \Big\{ \bar{x}_1^k < \bar{x}_2^k < \ldots < \bar{x}_{|\Mink|}^k :
				b(\bar{x}_j^k) = 0, \, - (k+1) b'(\bar{x}_j^k) \geq c(\bar{x}_j^k) \Big\}
			\subseteq \Minn.
		\end{gather*}
	\item By~\eqref{ieq:boundExpAway} and due to the structure of power layers
		the solution $u$ of problem~\eqref{prob:general} can be bounded by a
		constant away from the location of boundary and interior layers. Thus
		a layer-adaption of the mesh is only needed in the vicinity of 
		$\bar{x} \in \Mink \cup \Bpow \cup \Bexp$, i.e., in a $\delta$-neighborhood
		with $\delta > 0$ independent of $\eps$.
\end{itemize}

Now, for $\bar{x} \in \Bexp$ define a transition parameter $\tau_{\bar{x}}$ as given
in~\eqref{def:tau} with $\rho \geq k+1$. In the case that $\bar{x} \in \{\ua,\oa\} \setminus \Bexp$
we set $\tau_{\bar{x}} = 0$. We will assume in the following that
\begin{gather}
	\label{ass:tau}
	2 \tau_{\bar{x}} \leq \delta_{\bar{x},k}
		= \min\!\Big\{ |\bar{x}-y| : y \in \left(\{\ua,\oa\}\setminus \{\bar{x}\}\right)
			\cup \Mink \Big\}
	\qquad \text{for any } \bar{x} \in \Bexp
\end{gather}
which is the interesting case in practice, see also Remark~\ref{rem:equiMesh}.
In order to define a suitable mesh the following strategy could be pursued:
\begin{itemize}
	\item \textbf{Exponential boundary layer region:} If $\ua \in \Bexp$ define a fine
		mesh via a mesh generating function $\varphi$ as used for S-type meshes, see
		Section~\ref{subsec:S-typeMesh}, inside the interval $[\ua,\ua+\tau_\ua]$.
		Analogously, a fine mesh is define on $[\oa-\tau_\oa,\oa]$ if $\oa \in \Bexp$.
	\item \textbf{Power-type layer region:} In the vicinity of $\Bpow \cup \Mink$ build
		the mesh as given in Section~\ref{subsec:SunStynesMesh} with properly chosen
		parameter $\lambda$, i.e., $\lambda = 0$ for $\bar{x} \in \Bpow$ and
		$\lambda = \mu c(\bar{x})/|b'(\bar{x})|$, $\mu \in (0,1)$, for $\bar{x} \in \Mink$.
		Note that by construction we have $h_i^{-1} \snorm{b}_{\infty,[x_{i-1},x_i]} \leq C N$,
		see also Remark~\ref{rem:xihiInv}. Its appropriate to use this kind of mesh design on
		the whole interval $[\ua+\tau_\ua,\oa-\tau_\oa]$ if $\Bpow \cup \Mink \neq \emptyset$.
	\item \textbf{Rest of the interval:} In the rest of the interval the mesh should
		by chosen such that the mesh interval lengths satisfy $h_i \leq C N^{-1}$ and
		$h_i^{-1} \snorm{b}_{\infty,[x_{i-1},x_i]} \leq C N$.
\end{itemize}
The readers who would rather learn the mesh construction from specific examples are
referred to Section~\ref{sec:examples}. There we give specific layer-adapted meshes
for problems with different layers. Anyway these should clarify most of the
arising questions regarding the construction strategy.

\begin{remark}
	\label{rem:equiMesh}
	If $2 \tau_{\bar{x}} > \delta_{\bar{x},k}$ for any $\bar{x} \in \Bexp$,	we at
	worst gain that $\eps^{-1/2} \leq C \ln N$ which enables the use of an equidistant
	mesh. Indeed, even in the worst case a standard argumentation then yields the
	uniform error estimate
	\begin{gather*}
		\tnorm{u-u^N}_\eps \leq C \norm{u-u^I}_1 \leq  C N^{-k} \snorm{u^{(k+1)}}
			\leq C N^{-k} \left(1+\eps^{-(k+1)}\right) \leq C N^{-k} (\ln N)^{2(k+1)}.
	\end{gather*}
	on a mesh with $h_i \leq C N^{-1}$. In practice, where typically $0 < \eps \ll 1$, such a
	situation is uncommon since it would imply that $N$ is exponentially large
	compared to $\eps^{-1/2}$.
\end{remark}

\section{Error estimation}
\label{sec:errorEst}

We want to present a uniform error estimation for problem~\eqref{prob:general}
with~\eqref{prob:assNew} on a layer-adapted mesh as suggested in
Section~\ref{subsec:genMesh} now. For this recall that by Theorem~\ref{th:solDecomp1}
the solution of the problem can be split in an exponential boundary
layer part and a part that is smooth or at the most includes power-type
layers. Both solution parts will be estimated separately.

\begin{theorem}
	\label{th:errorGen}
	Let $u$ and $u^N \in V^N$ be the solutions of~\eqref{weakprob:general}
	and~\eqref{dprob:general}, respectively. Let $k \geq 1$ denote the
	ansatz order of the discrete space $V^N$ defined on a mesh as suggested
	in Section~\ref{subsec:genMesh}, in particular with $\rho \geq k+1$.
	Then, assuming~\eqref{ass:tau}, we have (in the worst case)
	\begin{gather*}
		\tnorm{u-u^N}_{\eps}
			\leq C \left((\KAK_0 + 1)N^{-1} + h + N^{-1} \max |\psi'| \right)^{k},
	\end{gather*}
	where $h := \max_{x_i \in (\ua,\ua+\tau_\ua] \cup (\oa-\tau_\oa,\oa]} h_i
		\leq C \sqrt{\eps} N^{-1} \max \varphi' \leq C \sqrt{\eps}$
	and $\KAK_0$ as in~\eqref{KAK} with $\lambda = 0$.
\end{theorem}

\begin{proof}
	The error estimation is based on the splitting~\eqref{splitError}, i.e.,
	$u-u^N = (u-u^I)+(u^I-u^N)$. Let $v^N = u^N -u^I$. In order to shorten
	the formulas we introduce the notation $I_\ua := (\ua,\ua+\tau_\ua)$ and
	$I_\oa := (\oa-\tau_\oa,\oa)$. From~\eqref{startEstimation},
	the solution decomposition, and	integration by parts we gain
	\begin{align*}
		\tnorm{u^N-u^I}_{\eps}^2
			& \leq B_\eps(u^N,v^N) - B_\eps(u^I,v^N) = B_\eps(u,v^N) - B_\eps(u^I,v^N) \\
			& = \eps \big((u-u^I)',(v^N)'\big) + \big(f(\cdot,u)-f(\cdot,u^I),v^N\big) + \big(b(S-S^I)',v^N\big) \\
			& \qquad - \big(b'(E-E^I),v^N\big)_{(\ua+\tau_\ua,\oa-\tau_\oa)}
				- \big(b(E-E^I),(v^N)'\big)_{(\ua+\tau_\ua,\oa-\tau_\oa)} \\
			& \qquad + \big(b(E_\ua-E_\ua^I)',v^N\big)_{(\oa-\tau_\oa,\oa)}
				+ \big(b(E_\oa-E_\oa^I)',v^N\big)_{(\ua,\ua+\tau_\ua)} \\
			& \qquad + \sum_{\bar{x} \in \{\ua,\oa\}} \Big\{
				(1-\mu_{\bar{x}}) \big(b(E_{\bar{x}}-E_{\bar{x}}^I)',v^N\big)_{I_{\bar{x}}} \\
			& \qquad \qquad - \mu_{\bar{x}} \left[ \big(b'(E_{\bar{x}}-E_{\bar{x}}^I),v^N\big)_{I_{\bar{x}}}
				+ \big(b(E_{\bar{x}}-E_{\bar{x}}^I),(v^N)'\big)_{I_{\bar{x}}} \right] \Big\},
	\end{align*}
	where $\mu_\ua, \mu_\oa \in [0,1]$. Using~\eqref{ieq:nonlinToLin} and applying
	Cauchy-Schwarz' inequality we obtain
	\begin{align*}
		\tnorm{u^N-u^I}_{\eps}^2
			& \leq C \underbrace{\left(\tnorm{S-S^I}_{\eps} + \norm{b(S-S^I)'}_0
				\right)}_{\mathrm{(i)}} \tnorm{v^N}_{\eps} \\
			& \qquad + C \underbrace{\left( \tnorm{E-E^I}_{\eps}
					+\big\|E_\ua-E_\ua^I\big\|_{0,(\ua,\ua+\tau_\ua)}
					+ \big\|E_\oa-E_\oa^I\big\|_{0,(\oa-\tau_\oa,\oa)}
				\right)}_{\mathrm{(ii)}} \norm{v^N}_{0} \\
			& \qquad + C \underbrace{\left( \big\|(E_\ua-E_\ua^I)'\big\|_{0,(\oa-\tau_\oa,\oa)}
					+ \big\|(E_\oa-E_\oa^I)'\big\|_{0,(\ua,\ua+\tau_\ua)}
				\right)}_{\mathrm{(iii)}} \norm{v^N}_{0} \\
			& \qquad + \underbrace{\big| \big(b(E-E^I),(v^N)'\big)_{(\ua+\tau_\ua,\oa-\tau_\oa)} \big|}_{
				\mathrm{(iv)}} \\
			& \qquad + \sum_{\bar{x} \in \{\ua,\oa\}} \underbrace{
				\min \Big\{\big| \big(b(E_{\bar{x}}-E_{\bar{x}}^I),(v^N)'\big)_{I_{\bar{x}}} \big| ,
				\big|\big(b(E_{\bar{x}}-E_{\bar{x}}^I)',v^N\big)_{I_{\bar{x}}} \big| \Big\}}_{
					\mathrm{(v)}}.
	\end{align*}
	
	The interpolation errors in $\mathrm{(i)}$ can be estimated by standard estimates
	if $\Bpow \cup \Mink= \emptyset$. Otherwise the mesh is properly adapted to the
	power-type layers in the vicinity of $\bar{x} \in \Bpow \cup \Mink$	and, thus,
	Lemma~\ref{le:IntErrorSunStynes} can be used there. Also note that the presence
	of a power-type layer at $\bar{x}$ correlates to $|b(x)| \leq C |x-\bar{x}|$
	and, furthermore, that $S$ and its derivatives can be bounded by a constant independent
	of $\eps$ in $[\ua,\ua+\tau_\ua] \cup [\oa-\tau_\oa,\oa]$.
	
	In order to bound $\mathrm{(ii)}$ we simply have to apply Lemma~\ref{le:IntErrorS-type}.
	Here recall that the exact mesh structure in $[\ua+\tau_\ua,\oa-\tau_\oa]$ is not important,
	see also Remark~\ref{rem:coarseMeshPart}. The term $\mathrm{(iii)}$	can be treated by
	standard interpolation theory since~\eqref{ieq:boundExpAway} implies that $E_\ua$ and
	$E_\oa$ and their derivatives can be $\eps$-uniformly bounded in $[\oa,\oa-\tau_\oa]$
	and $[\ua,\ua+\tau_\ua]$, respectively.
	
	The mesh construction guarantees that $h_i^{-1} \snorm{b}_{\infty,[x_{i-1},x_i]} \leq C N$
	when $x_i \in (\ua+\tau_{\ua},\oa-\tau_{\oa}]$. Therefore, combining H\"{o}lder's inequality,
	an inverse inequality~\eqref{invIneq}, and Lemma~\ref{le:IntErrorS-type} we obtain
	\begin{align*}
		\mathrm{(iv)}
			& = \big| \big(b(E-E^I),(v^N)'\big)_{(\ua+\tau_\ua,\oa-\tau_\oa)} \big| \\
			& \leq C \norm{E-E^I}_{\infty,[\ua+\tau_\ua,\oa-\tau_\oa]}
				\sum_{x_i \in (\ua+\tau_\ua,\oa-\tau_\oa]} \int_{x_{i-1}}^{x_i} \snorm{b} \snorm{(v^N)'} \\
			& \leq C \norm{E-E^I}_{\infty,[\ua+\tau_\ua,\oa-\tau_\oa]}
				\sum_{x_i \in (\ua+\tau_\ua,\oa-\tau_\oa]}
					\snorm{b}_{\infty,[x_{i-1},x_i]} h_i^{-1} \int_{x_{i-1}}^{x_i} \snorm{v^N} \\
			& \leq C \norm{E-E^I}_{\infty,[\ua+\tau_\ua,\oa-\tau_\oa]}
				N \norm{v^N}_{0,(\ua+\tau_\ua,\oa-\tau_\oa)} \\
			& \leq C N^{-\rho+1} \tnorm{v^N}_{\eps}.
	\end{align*}
	
	It remains to study $\mathrm{(v)}$. There are two cases.
	\begin{enumerate}[(I)]
	\item If there is an
		exponential boundary layer of width $\KAO(\eps \log 1/\eps)$ located at $\bar{x}$ (first
		case in~\eqref{def:Ephi1}), then Cauchy-Schwarz' inequality and Lemma~\ref{le:IntErrorS-type}
		yield
		\begin{align*}
			\mathrm{(v)} \leq \big| \big(b(E_{\bar{x}}-E_{\bar{x}}^I),(v^N)'\big)_{I_{\bar{x}}} \big|
				& \leq C \norm{E_{\bar{x}}-E_{\bar{x}}^I}_{0,I_{\bar{x}}} \snorm{v^N}_{1} \\
				& \leq C \eps^{1/2} \left(N^{-1} \max |\psi'|\right)^{k} \snorm{v^N}_{1} \\
				& \leq C \left(N^{-1} \max |\psi'|\right)^{k} \tnorm{v^N}_\eps
		\end{align*}
		where we assumed that $N^{-1} \max|\psi'| (\ln N)^{1/2} \leq C$ or alternatively
		$\rho \geq k +3/2$.
	\item Otherwise there is an exponential boundary layer of width
		$\KAO(\sqrt{\eps} \log 1/\sqrt{\eps})$ located at $\bar{x}$ (second case
		in~\eqref{def:Ephi1}). But then we have $b(\bar{x}) = b'(\bar{x}) = 0$
		and, thus, $|b(x)| \leq C |x-\bar{x}|^2$. By H\"{o}lder's inequality
		and Lemma~\ref{le:xTimesIntErrorS-type} we obtain
		\begin{align*}
		\mathrm{(v)} \leq \left|\big(b(E_{\bar{x}}-E_{\bar{x}}^I)',v^N\big)\right|_{I_{\bar{x}}}
			& \leq C \norm{|x-\bar{x}|^2 (E_{\bar{x}}-E_{\bar{x}}^I)'}_{\infty,I_{\bar{x}}}
				\left(\meas(I_{\bar{x}})\right)^{1/2} \norm{v^N}_{0} \\
			& \leq C (4\eps)^{1/2} \left(N^{-1} \max |\psi'|\right)^{k}
				\underbrace{\tau_{\bar{x}}^{1/2}}_{\leq C} \norm{v^N}_{0}.
		\end{align*}
		Note that an argumentation as in (I) would not work here since from the
		length of the interval, i.e., from $\left(\meas(I_{\bar{x}})\right)^{1/2} = \tau_{\bar{x}}^{1/2}
			\leq C \eps^{1/4} (\ln N)^{1/2}$, we would get a factor $\eps^{1/4}$ only.
	\end{enumerate}
	Combining all the above estimates we easily complete the proof.
\end{proof}

\begin{remark}
	Note that from~\eqref{ieq:KAKtoLN} we have $\KAK_0+1 \leq C \ln N$.
	Moreover, if we choose a Shishkin mesh in the exponential boundary
	layer region we have $\max \varphi' + \max |\psi'| \leq C \ln N$. So,
	in the setting of Theorem~\ref{th:errorGen} we can always guarantee
	that
	\begin{gather*}
		\tnorm{u-u^N}_{\eps}
			\leq C \left(N^{-1} \ln N \right)^{k}.
	\end{gather*}
\end{remark}

\section{Some examples of linear problems with different layers}
\label{sec:examples}

In this section we want to study some linear problems whose solutions
exhibit different types of layers. The problems of type~\eqref{prob:generalLin}
are stated in such a way that~\eqref{prob:assOld} is satisfied.
We will assume (without loss of generality for sufficiently small $\eps$,
see Lemma~\ref{le:transCmBsH}) that
\begin{gather}
	\label{prob:standAss}
	c(x) \geq 0, \qquad c(x) - \tfrac{1}{2} b'(x) \geq \tilde{\gamma} > 0.
\end{gather}
Note that this assumption will not be stated explicitly in the problem
descriptions.

For each problem we give a layer-adapted mesh constructed as suggested
in Section~\ref{subsec:genMesh} and state an error estimate. The error
bounds are slightly more precise than the one of Theorem~\ref{th:errorGen}
since we know more about the structure of the problems. But anyway a similar
proof could be used.

\subsection{Problem with a repulsive boundary turning point}

Consider the singularly perturbed boundary value problem~\eqref{prob:generalLin}
on $\oI := [0,1]$ with
\begin{gather}
	\label{prob:RepBouTPP}
	b(x) = x a(x), \qquad a(x) \geq \alpha > 0, \qquad c(0) > 0, \qquad c(0) - \tfrac{1}{2} b'(0) > 0.
\end{gather}
Obviously, the problem has a repulsive boundary turning point at $x = 0$.

From Section~\ref{sec:AprioriAndSolDecomp} we know that the solution $u$ of~\eqref{prob:generalLin}
with~\eqref{prob:standAss}, \eqref{prob:RepBouTPP} can be bounded as
\begin{gather*}
	\snorm{u^{(l)}(x)}
		\leq C \left(1 + \eps^{\lambda/2}\left(\eps^{1/2}+x\right)^{-\lambda-l}
			+ \eps^{-l} e^{-b(1)(1-x)/\eps} \right)
\end{gather*}
with $0 \leq \lambda < c(0)/|b'(0)| = c(0)/a(0)$. Moreover, it can be decomposed
such that $u = S + E$ where for $q \in \NN$ and $l = 0, \ldots, q$ we have
\begin{gather*}
	\snorm{S^{(l)}(x)}
		\leq C \left(1 + \eps^{\lambda/2}\left(\eps^{1/2}+x\right)^{-\lambda-l} \right)
		\leq C \left(1 + \left(\eps^{1/2}+x\right)^{-l} \right)
\end{gather*}
and with $\beta = b(1) = a(1) > 0$
\begin{gather*}
	\snorm{E^{(l)}(x)}
		\leq C \eps^{-l} e^{-\beta(1-x)/\eps}.
\end{gather*}

We next present a possible layer-adapted mesh. In order to treat both the
exponential boundary layer at $x=1$ and the power-type boundary layer at $x=0$,
the two mesh types presented in Section~\ref{sec:mesh} are combined. A transition
parameter $\tau$ is defined by
\begin{gather*}
	\tau = \rho \frac{\eps}{\beta} \ln N
\end{gather*}
where $\rho > 0$. We will assume that $\tau \leq 1/2$, that is~\eqref{ass:tau} holds true.

Let $\varphi$ denote a mesh generating function as used in Section~\ref{subsec:S-typeMesh}
and $\psi$ its associated mesh characterizing function. Furthermore, let $\Phi_{N/2,0}$
be the mesh generating function as introduced in Section~\ref{subsec:SunStynesMesh}.
Then the mesh points are defined by
\begin{gather}
	\label{mesh:RepBouTpp}
	x_i =
	\begin{cases}
		(1-\tau) \Phi_{N/2,0}(t_i),
			& t_i = 2i/N, \, i = 0,\ldots, N/2, \\
		1-\rho\tfrac{\eps}{\beta} \varphi(t_i),
			& t_i = (N-i)/N, \, i = N/2,\ldots,N.
	\end{cases}
\end{gather}
Note that the mesh is somewhat similar to an S-type mesh. Of course in $[1-\tau,1]$ we have
the fine part. However, the ``coarse'' part in $[0,1-\tau]$ is not simply (quasi) uniform
but also adapted to the power-type layer.

\begin{corollary}
	\label{cor:ErrorRepBouTPP}
	Under the assumptions of Theorem~\ref{th:errorGen}, in particular $\rho \geq k+1$, we have
	on the mesh~\eqref{mesh:RepBouTpp}
	\begin{gather*}
		\tnorm{u-u^N}_{\eps}
			\leq C \left((\KAK + 1)N^{-1} + h + N^{-1} \max |\psi'| \right)^{k},
	\end{gather*}
	where $h := \max_{N/2 +1 \leq i \leq N} h_i \leq C \eps N^{-1} \max\varphi' \leq C \eps$.
\end{corollary}

\subsection{Problem with an attractive and a multiple boundary turning point}

Consider the singularly perturbed boundary value problem~\eqref{prob:generalLin}
on $\oI := [0,1]$ with
\begin{gather}
	\label{prob:AttAndMultBouTPP}
	\begin{aligned}
		& b(x) = -x(1-x)^\ell a(x),\, \ell\geq 2, \quad
		& a(x) \geq \alpha & > 0, \\
		& c(0) - \tfrac{1}{2} b'(0) \geq c(0) > 0, \quad
		& c(1) - \tfrac{1}{2} b'(1) = c(1) & > 0.
	\end{aligned}
\end{gather}
The problem has an attractive boundary turning point at $x = 0$ and
an multiple boundary turning point at $x=1$.

From Section~\ref{sec:AprioriAndSolDecomp} we know that the solution $u$ of~\eqref{prob:generalLin}
with~\eqref{prob:standAss}, \eqref{prob:AttAndMultBouTPP} can be bounded as
\begin{gather*}
	\snorm{u^{(l)}(x)}
		\leq C \left(1 + \left(\eps^{1/2}+x\right)^{\lambda-l} + \eps \left(\eps^{1/2}+x\right)^{-l-2}
			+ \eps^{-l/2} e^{-\sqrt{c(1)}(1-x)/\sqrt{\eps}} \right)
\end{gather*}
with $0 \leq \lambda < c(0)/|b'(0)| = c(0)/a(0)$. Moreover, it can be decomposed
such that $u = S + E$ where for $q \in \NN$ and $l = 0, \ldots, q$ we have
\begin{gather*}
	\snorm{S^{(l)}(x)}
		\leq C \left(1 + \left(\eps^{1/2}+x\right)^{\lambda-l} + \eps \left(\eps^{1/2}+x\right)^{-l-2} \right)
		\leq C \left(1 + \left(\eps^{1/2}+x\right)^{-l} \right)
\end{gather*}
and with $\beta = \sqrt{c(1)} > 0$
\begin{gather*}
	\snorm{E^{(l)}(x)}
		\leq C \eps^{-l/2} e^{-\beta(1-x)/\sqrt{\eps}}.
\end{gather*}

A layer-adapted mesh can be constructed as follows. Define
a transition parameter $\tau$ by
\begin{gather*}
	\tau = \rho \frac{\sqrt{\eps}}{\beta} \ln N
\end{gather*}
where $\rho > 0$. Also here we shall suppose that $\tau \leq 1/2$ in order to fulfill~\eqref{ass:tau}.

Similar to the mesh of the last subsection the mesh points now are defined by
\begin{gather}
	\label{mesh:AttAndMultBouTpp}
	x_i =
	\begin{cases}
		(1-\tau) \Phi_{N/2,0}(t_i),
			& t_i = 2i/N, \, i = 0,\ldots, N/2, \\
		1-\rho\tfrac{\sqrt{\eps}}{\beta} \varphi(t_i),
			& t_i = (N-i)/N, \, i = N/2,\ldots,N.
	\end{cases}
\end{gather}

\begin{corollary}
	Under the assumptions of Theorem~\ref{th:errorGen}, in particular $\rho \geq k+1$, we have
	on the mesh~\eqref{mesh:AttAndMultBouTpp}
	\begin{gather*}
		\tnorm{u-u^N}_{\eps}
			\leq C \left((\KAK + 1)N^{-1} + h \right)^k + C \eps^{1/4} \left(N^{-1} \max |\psi'| \right)^{k},
	\end{gather*}
	where $h := \max_{N/2 +1 \leq i \leq N} h_i \leq C \sqrt{\eps} N^{-1} \max\varphi' \leq C \sqrt{\eps}$.
\end{corollary}

\subsection{Problem with interior and boundary turning points}
Consider the singularly perturbed boundary value problem~\eqref{prob:generalLin}
on $\oI := [-1,1]$ with
\begin{gather}
	\label{prob:IntAndBouTPP}
	\begin{aligned}
		b(x) & = -(x+1)x(x-\tfrac{1}{2})(x-\tfrac{27}{30})^3 a(x), 
		& a(x) & \geq \alpha > 0,  \\
		c(x) & \geq \gamma > 0, 
		& c(x) - \tfrac{1}{2} b'(x) & \geq \tilde{\gamma} > 0.
	\end{aligned}
\end{gather}
The problem has a repulsive boundary turning point at $x = -1$, and
some interior turning points at $x=0$, $\frac{1}{2}$, and $\frac{27}{30}$.
An interior layer only occurs at $x = 0$. At $x=-1$ there is a power-type
boundary layer and at $x=1$ there is actually no (boundary) layer.

From Section~\ref{sec:AprioriAndSolDecomp} we know that the solution $u$ of~\eqref{prob:generalLin}
with~\eqref{prob:IntAndBouTPP} can be bounded as
\begin{align*}
	\snorm{u^{(l)}(x)}
		& \leq C \left(1 + \eps^{\lambda/2}\left(\eps^{1/2}+(x+1)\right)^{-\lambda-l}
			+ \left( \eps^{1/2} + |x| \right)^{\lambda_1-l} \right) \\
		&  \leq C \left(1 + \left(\eps^{1/2}+(x+1)\right)^{-l}
			+ \left( \eps^{1/2} + |x| \right)^{\lambda_1-l} \right)
\end{align*}
with $0 \leq \lambda < c(-1)/|b'(-1)|$ and $0 \leq \lambda_1 < c(0)/|b'(0)|$.

A suitable mesh could be defined in the following way
\begin{gather}
	\label{mesh:IntAndBouTPP}
	x_i =
	\begin{cases}
		-1 + \frac{1}{2} \Phi_{N/4,0}(t_i),
			& t_i = 4i/N, \, i = 0,\ldots, N/4, \\
		- \frac{1}{2} \Phi_{N/4,\lambda_1}(t_i),
			& t_i = 2(N-2i)/N, \, i = N/4,\ldots, N/2, \\
		\Phi_{N/2,\lambda_1}(t_i),
			& t_i = (2i-N)/N, \, i = N/2,\ldots,N,
	\end{cases}
\end{gather}
where we supposed that the ``difficult'' case $\lambda_1 \in [0,k+1)$
is present.

\begin{corollary}
	\label{cor:ErrorIntAndBouTPP}
	Under the assumptions of Theorem~\ref{th:errorGen}, we have
	on the mesh~\eqref{mesh:IntAndBouTPP}
	\begin{gather*}
		\tnorm{u-u^N}_{\eps}
			\leq C \left((\KAK_0 + 1)N^{-1}\right)^{k}.
	\end{gather*}
\end{corollary}

\subsection{Problem with two exponential boundary layers of different width}
Consider the singularly perturbed boundary value problem~\eqref{prob:generalLin}
on $\oI := [0,1]$ with
\begin{gather}
	\label{prob:TwoExpLayerTPP}
	b(x) = x^\ell a(x), \, \ell \geq 2 \qquad
	a(x) \geq \alpha > 0, \qquad
	c(0) - \tfrac{1}{2} b'(0) = c(0) > 0.
\end{gather}
The problem has an exponential boundary layer of width $\KAO(\sqrt{\eps}\log 1/\sqrt{\eps})$
at the multiple boundary turning point at $x = 0$ and an exponential boundary
layer of width $\KAO(\eps\log 1/\eps)$ at the outflow boundary $x = 1$.

From Section~\ref{sec:AprioriAndSolDecomp} we know that the solution $u$ of~\eqref{prob:generalLin}
with~\eqref{prob:standAss}, \eqref{prob:TwoExpLayerTPP} can be bounded as
\begin{gather*}
	\snorm{u^{(l)}(x)}
		\leq C \left(1 + \eps^{-l/2} e^{-\sqrt{c(0)}x/\sqrt{\eps}}
			+ \eps^{-l}e^{-b(1)(1-x)/\eps} \right)\!.
\end{gather*}
Moreover, it can be decomposed such that $u = S + E_0 + E_1$ where for $q \in \NN$
and $l = 0, \ldots, q$ we have
\begin{gather*}
	\snorm{S^{(l)}(x)} \leq C, \qquad
	\snorm{E_0^{(l)}(x)}
		\leq C \eps^{-l/2} e^{-\beta_0 x/\sqrt{\eps}}, \qquad
	\snorm{E_1^{(l)}(x)}
		\leq C \eps^{-l} e^{-\beta_1(1-x)/\eps}
\end{gather*}
with $\beta_0 = \sqrt{c(0)} > 0$ and $\beta_1 = b(1)= a(1) > 0$.

We use a layer-adapted S-type mesh. The transition parameters are
defined by
\begin{gather*}
	\tau_0 = \rho \frac{\sqrt{\eps}}{\beta_0} \ln N, \qquad
	\tau_1 = \rho \frac{\eps}{\beta_1} \ln N.
\end{gather*}
We shall assume that $\max\{\tau_0, \tau_1\} \leq \frac{1}{4}$ in the following
which guarantees~\eqref{ass:tau}. A suitable mesh is given by the mesh points
\begin{gather}
	\label{mesh:TwoExpLayerTPP}
	x_i =
	\begin{cases}
		\rho\tfrac{\sqrt{\eps}}{\beta_0} \varphi(t_i),
			& t_i = 2i/N, \, i = 0,\ldots,N/4, \\
		\tau_0 + (1-\tau_0 - \tau_1) t_i ,
			& t_i = (4i-N)/(2N), \, i = N/4,\ldots, 3N/4, \\
		1-\rho\tfrac{\eps}{\beta_1} \varphi(t_i),
			& t_i = 2(N-i)/N, \, i = 3N/4,\ldots,N.
	\end{cases}
\end{gather}

\begin{corollary}
	\label{cor:ErrorTwoExpLayerTPP}
	Under the assumptions of Theorem~\ref{th:errorGen},	in particular $\rho \geq k+1$,
	we have on the mesh~\eqref{mesh:TwoExpLayerTPP}
	\begin{gather*}
		\tnorm{u-u^N}_{\eps}
			\leq C \left(N^{-1} + h + N^{-1} \max |\psi'| \right)^{k},
	\end{gather*}
	where $h := \max_{1 \leq i \leq N/4 \, \vee\, 3N/4+1 \leq i \leq N} h_i \leq C \sqrt{\eps} N^{-1} \max\varphi' \leq C \sqrt{\eps}$.
\end{corollary}

\begin{appendix}
	\section{Problem transformation}
	\label{app:probTrans}
	In this section transformations of the problem
	\begin{gather}
	\tag{\ref{prob:generalLin}}
	\begin{aligned}
		Lu := - \eps u''(x) + b(x)u'(x) + c(x)u(x) & = f(x), \qquad \text{for } x \in \oI := [\ua,\oa], \\
		u(\ua)= \nu_{-}, \quad u(\oa) & = \nu_+,
	\end{aligned}
	\end{gather}
	with $\eps \in (0,\eps_0]$, $\eps_0 > 0$ sufficiently small, problem data $b,c,f$ sufficiently smooth, and
	\begin{align}
		\text{for all } x \in \oI: \qquad b(x) = 0 \quad & \Longrightarrow \quad
			c(x)>0, \tag{\ref{prob:ass1}} \\
		\text{for all } x \in \oI: \qquad b(x) = 0 \quad & \Longrightarrow \quad
			\!\left(c - \tfrac{1}{2}b'\right)(x) >0 \tag{\ref{prob:ass2}}
	\end{align}
	shall be studied. These transformations are used to guarantee certain
	properties of the problem data. 
	
	We define the compact set of all zeros of $b$ by
	$M:= \{x \in \oI : b(x) = 0 \} = \oI \cap b^{-1}(\{0\})$.
	Applying the extreme value theorem we conclude that
	\begin{gather}
		\label{ieq:c0bound}
		\min_{x \in M} c(x) \geq c_0 > 0
		\quad \text{or} \quad
		\min_{x \in M} \left(c-\tfrac{1}{2}b'\right)(x) \geq \tilde{c}_0 > 0,
		\quad \text{respectively}.
	\end{gather}
	Moreover, there is a $\delta_0 > 0$ such that
	\begin{gather}
		\label{ieq:delta0}
		\min_{x \in B(M,\delta_0)} c(x) \geq \frac{c_0}{2} > 0
		\quad \text{or} \quad
		\min_{x \in B(M,\delta_0)} \left(c-\tfrac{1}{2}b'\right)(x) \geq \frac{\tilde{c}_0}{2} > 0,
		\quad \text{respectively}.
	\end{gather}
	Note that $M$ is the union of a finite number of closed disjoint subintervals of $\oI$ only.
	
	\begin{lemma}
		\label{le:pDef}
		Let $\delta > 0$. There is a function $p$ such that
		\begin{itemize}
			\item $p \in C^\infty(\oI)$,
			\item $p \geq 1$,
			\item $0 \leq \sgn(b(x)) p'(x) \leq 1$ \quad and \quad
				$p'(x)= \begin{cases}
							0			& \text{for } x \in \oI \cap B(M,\delta/3), \\
							\sgn(b(x))	& \text{for } x \in \oI \setminus B(M,\delta).
						\end{cases}$
		\end{itemize}
		Here $\sgn(x)$ denotes the sign of $x$.
	\end{lemma}
	\begin{proof}
		Define a continuous, piecewise linear function $\hat{p}$ such that
		\begin{gather*}
			\hat{p}'(x) = \sgn(b(x)) \one_{\oI \setminus B(M,\frac{2}{3}\delta)}
			\qquad \text{for almost all $x \in [\ua,\oa]$,}
		\end{gather*}
		where $\one_A$ denotes the indicator function of $A \subseteq \RR$. Then
		extent this function continuously to $[\ua-1,\oa+1]$ such that
		\begin{gather*}
			\hat{p}'(x) = \lim_{y \searrow \ua} \hat{p}'(y)
			\quad x \in [\ua-1,\ua)
			\qquad \text{and} \qquad
			\hat{p}'(x) = \lim_{y \nearrow \oa} \hat{p}'(y)
			\quad x \in (\oa,\oa+1].
		\end{gather*}
		The wanted function $p$ can be easily constructed with the help of
		a convolution of $\hat{p}$ with a suitable mollifier and if necessary by adding
		a suitable constant.
	\end{proof}
	
	\begin{lemma}
		\label{le:transC}
		Consider problem~\eqref{prob:generalLin} with~\eqref{prob:ass1} and $c(x) \geq 0$.
		Then we may assume without loss of generality that there is a constant $\gamma >0$
		such that $c(x) \geq \gamma > 0$ for all $x \in \oI$.
	\end{lemma}
	\begin{proof}
		Using the auxiliary function $p$ with $\delta = \delta_0$,
		we define a new operator
		\begin{gather*}
			\tilde{L}_\kappa v := - \eps v'' + \tilde{b}_\kappa(x)v' + \tilde{c}_\kappa(x) v \qquad (\kappa > 0)
		\end{gather*}
		where
		\begin{gather*}
			\tilde{b}_\kappa(x):=b(x)-2 \eps \kappa p'(x),
			\qquad \tilde{c}_\kappa(x) := c(x)-\eps \kappa p''(x) + b(x) \kappa p'(x) - \eps \kappa^2 (p'(x))^2.
		\end{gather*}
		Noting that
		\begin{align*}
			\left(e^{-\kappa p(x)}v(x)\right)'
				& = -\kappa p'(x) e^{-\kappa p(x)}v(x) + e^{-\kappa p(x)} v'(x), \\
			\left(e^{-\kappa p(x)}v(x)\right)''
				& = \left(-\kappa p''(x) + \kappa^2 (p'(x))^2 \right) e^{-\kappa p(x)}v(x)
					- 2\kappa p'(x) e^{-\kappa p(x)}v'(x) + e^{-\kappa p(x)} v''(x),
		\end{align*}
		it is easy to verify that
		\begin{align*}
			\tilde{L}_\kappa \left(e^{-\kappa p(x)} v(x)\right)
				& = e^{-\kappa p(x)} \Big[-\eps\left(\left(-\kappa p''(x)+\kappa^2(p'(x))^2\right)v(x)-2\kappa p'(x)v'(x)+v''(x)\right) \\
					& \qquad + \left(-b(x)\kappa p'(x)+2\eps\kappa^2 (p'(x))^2\right)v(x)+\left(b(x)-2\eps\kappa p'(x)\right)v'(x) \\
					& \qquad + \left(c(x)-\eps \kappa p''(x) + b(x) \kappa p'(x) - \eps \kappa^2 (p'(x))^2\right) v(x) \Big] \\
				& = e^{-\kappa p(x)} L v(x).
		\end{align*}
		We distinguish three cases:
		\begin{enumerate}
			\item[(i)] $x \in \oI \cap B(M,\delta_0/3)$: \\
				By Lemma~\ref{le:pDef} and~\eqref{ieq:c0bound} we have $\tilde{b}_\kappa(x) = b(x)$ and
				$\tilde{c}_\kappa(x) = c(x) \geq c_0/2 > 0$.
			\item[(ii)] $x \in \oI \cap (B(M,\delta_0) \setminus B(M,\delta_0/3))$: \\
				Since $\min \{|b(x)| : x \in \oI \setminus B(M,\delta_0/3) \} =: b_0 > 0$,
				we can ensure that $\sgn(b(x)) = \sgn(\tilde{b}_\kappa(x)) \neq 0$ by choosing $\kappa>0$
				sufficiently small ($0 < \kappa \leq b_0/(4 \eps)$). Furthermore, by Lemma \ref{le:pDef} we have
				\begin{align*}
					\tilde{c}_\kappa(x)
						\geq c(x) - \kappa \eps \left(p''(x) + \kappa p'(x)^2 \right)
						\geq c(x) - \kappa \eps \left(\max |p''| + \kappa \right) \geq c_0/4 > 0
				\end{align*}
				for $0< \kappa \leq \min\left\{1, c_0/(4\eps) \left(1 + \max |p''|\right)^{-1}\right\}$.
			\item[(iii)] $x \in \oI \setminus B(M,\delta_0)$: \\
				Again, because of $\min \{|b(x)| : x \in \oI \setminus B(M,\delta_0/3) \} =: b_0 > 0$,
				we can ensure that $\sgn(b(x)) = \sgn(\tilde{b}_\kappa(x)) \neq 0$ by choosing $\kappa$
				sufficiently small ($0 < \kappa \leq b_0/(4 \eps)$). 
				Exploiting the properties of $p$ (cf.~Lemma~\ref{le:pDef}) we obtain
				for $0 < \kappa \leq b_0/(2\eps)$ that
				\begin{gather*}
					\tilde{c}_{\kappa}(x) = c(x) + \kappa \left(|b(x)| - \eps \kappa\right)
						\geq c(x) + \kappa \left(b_0 - \eps \kappa\right) \geq \kappa b_0/2.
				\end{gather*}
		\end{enumerate}
		Summarizing, for $0 < \kappa \leq \min\{1,b_0/(4\eps),c_0/(4\eps) (1+\max |p''|)^{-1}\}$
		the statement is proven with $\gamma = \min \{c_0/4, \kappa b_0/2\} > 0$.
	\end{proof}
	
	\begin{lemma}
		\label{le:transCmBsH}
		Consider problem~\eqref{prob:generalLin} with~\eqref{prob:assOld}. Then for sufficiently
		small $\eps_0 > 0$ we may assume without loss of generality that there are constants
		$\gamma, \tilde{\gamma} >0$ such that $c(x) \geq \gamma >0$ and 
		$\left(c-\tfrac{1}{2}b'\right)(x) \geq \tilde{\gamma} > 0$ for all $x \in \oI$.
	\end{lemma}
	\begin{proof}
		Let $\delta_0$ be chosen such that both inequalities in~\eqref{ieq:delta0} hold.
		As in the proof of Lemma~\ref{le:transC} we use $p$ with $\delta = \delta_0$ to define
		the modified operator
		\begin{gather*}
			\tilde{L}_\kappa v := - \eps v'' + \tilde{b}_\kappa(x)v' + \tilde{c}_\kappa(x) v \qquad (\kappa > 0)
		\end{gather*}
		where
		\begin{gather*}
			\tilde{b}_\kappa(x):=b(x)-2 \eps \kappa p'(x),
			\qquad \tilde{c}_\kappa(x) := c(x)-\eps \kappa p''(x) + b(x) \kappa p'(x) - \eps \kappa^2 (p'(x))^2.
		\end{gather*}
		We already had seen that
		\begin{gather*}
			\tilde{L}_\kappa \left(e^{-\kappa p(x)} v(x)\right)
				= e^{-\kappa p(x)} L v(x).
		\end{gather*}
		
		Now, we address the lower bound for $(c-\frac{1}{2}b')(x)$. 
		Set $\tilde{\gamma}_0 := \min_{x \in \oI} (c-\frac{1}{2}b')(x)$.
		We can assume that $\tilde{\gamma}_0 \leq 0$. Additionally, set
		$\gamma_0 := \min_{x \in \oI} c(x)$. We may also allow $\gamma_0 \leq 0$
		in the following.
		
		Three different cases have to be considered:
		\begin{enumerate}
			\item[(i)] $x \in \oI \cap B(M,\delta_0/3)$: \\
				By Lemma~\ref{le:pDef}, \eqref{ieq:c0bound}, and~\eqref{ieq:delta0} we have $\tilde{b}_\kappa(x) = b(x)$
				and $\tilde{c}_\kappa(x) = c(x) \geq c_0/2 > 0$. Furthermore, we obtain
				$\left(\tilde{c}_\kappa - \tfrac{1}{2}b_\kappa'\right)\!(x)
					= \left(c-\tfrac{1}{2}b'\right)\!(x) \geq \frac{\tilde{c}_0}{2} > 0$.
			\item[(ii)] $x \in \oI \cap (B(M,\delta_0) \setminus B(M,\delta_0/3))$: \\
				Since $\min \{|b(x)| : x \in \oI \setminus B(M,\delta_0/3) \} =: b_0 > 0$ we
				obtain with Lemma \ref{le:pDef} that
				\begin{align*}
					|\tilde{b}_\kappa(x)| \geq |b(x)|- 2 \eps \kappa \geq b_0/2 > 0,
					\qquad \sgn(\tilde{b}_\kappa(x)) = \sgn(b(x)) \neq 0
				\end{align*}
				for $0 < \eps \leq b_0/(4\kappa)$. Likewise we get
				\begin{align*}
					\tilde{c}_\kappa(x)
						\geq c(x) - \kappa \eps \left(p''(x) + \kappa p'(x)^2 \right)
						\geq c(x) - \kappa \eps \left(\max |p''| + \kappa \right) \geq c_0/4 > 0
				\end{align*}
				for $0< \eps \leq c_0/(4\kappa (\max |p''|+\kappa))$ and
				\begin{gather*}
					(\tilde{c}_\kappa -\tfrac{1}{2} \tilde{b}'_\kappa)(x)
						\geq (c- \tfrac{1}{2}b')(x) - \eps \kappa^2
						\geq \tilde{c}_0/4 > 0
				\end{gather*}
				for $0 < \eps \leq \tilde{c}_0/(4\kappa^2)$.
			\item[(iii)] $x \in \oI \setminus B(M,\delta_0)$: \\
				As in (ii) we can ensure that $\sgn(b(x)) = \sgn(\tilde{b}_\kappa(x)) \neq 0$
				as well as $|\tilde{b}_\kappa(x)| \geq b_0/2 > 0$ for $\eps$ sufficiently small
				($0 < \eps \leq b_0/(4 \kappa)$). By the properties of $p$ (cf.~Lemma~\ref{le:pDef})
				we furthermore obtain for $0 < \eps \leq 1/(\kappa^2)$.
				\begin{align*}
					\tilde{c}_\kappa(x)
					  & = c(x) +  \left(\kappa |b(x)| - \eps \kappa^2 \right)
						\geq \gamma_0 + \kappa b_0 - \eps \kappa^2
						\geq (\gamma_0 - 1) +\kappa b_0
				\intertext{and}
					(\tilde{c}_\kappa -\tfrac{1}{2} \tilde{b}'_\kappa)(x)
					  & = (c- \tfrac{1}{2}b')(x) + \kappa |b(x)| - \eps \kappa^2
						\geq \tilde{\gamma}_0 + \kappa b_0 - \eps \kappa^2
						\geq (\tilde{\gamma}_0 - 1) + \kappa b_0.
				\end{align*}
				Hence, choosing $\kappa \geq (\max\{0,1-\gamma_0,1-\tilde{\gamma}_0\}+\max\{c_0,\tilde{c}_0\}/4)/b_0 > 0$
				yields
				\begin{gather*}
					\tilde{c}_\kappa(x) \geq c_0/4 > 0,
					\qquad
					(\tilde{c}_\kappa -\tfrac{1}{2} \tilde{b}'_\kappa)(x) \geq \tilde{c}_0/4 > 0.
				\end{gather*}
		\end{enumerate}
		Summarizing, choosing $\kappa \geq (\max\{0,1-\gamma_0,1-\tilde{\gamma}_0\}+\max\{c_0,\tilde{c}_0\}/4)/b_0 > 0$
		the statement is proven with $\gamma = c_0/4$ and $\tilde{\gamma} = \tilde{c}_0/4$ for all
		$0 < \eps \leq \eps_0$ with $\eps_0$ sufficiently small. More detailed, it would be sufficient if
		\begin{gather*}
			0 < \eps_0 \leq \min\{1/\kappa^2, b_0/(4 \kappa),c_0/(4\kappa(\max|p''|+\kappa)),\tilde{c}_0/(4\kappa^2)\}.
		\end{gather*}
	\end{proof}
	
	\section{Alternative solution decomposition}
	\label{app:solDecomp}
	
	In this section we want to present an alternative way to derive a
	decomposition of the solution of~\eqref{prob:generalLin} with~\eqref{prob:assNew}.
	We first state the final result.
	
	\begin{theorem}[Solution decomposition]
		\label{th:solDecomp2}
		Let $q \in \NN$ and suppose that $b \in C^{2q,1}(\oI)$, $c \in C^{2q-1,1}(\oI)$, and $f \in C^q(\oI)$.
		Then the solution $u$ of~\eqref{prob:generalLin} with~\eqref{prob:assNew} has the representation
		\begin{gather*}
			u = S + E_\ua + E_\oa + \eps^q R,
		\end{gather*}
		where for all $k = 0, \ldots ,q$ and $x \in \oI$
		\begin{gather}
		\label{ieq:Sbound2}
		\begin{aligned}
			|S^{(k)}(x)|
				& \leq C \bigg( 1 + \phi^S_\ua\big((x-\ua),k,b(\ua),b'(\ua),\eps\big)
					+ \phi^S_\oa\big((\oa-x),k,-b(\oa),b'(\oa),\eps\big) \\
				& \hspace{18em} + \sum_j \left( \eps^{1/2} + |x-\bar{x}_j| \right)^{\lambda_j-k} \bigg),
		\end{aligned}
		\end{gather}
		with $0 < \lambda_j < c(\bar{x}_j)/|b'(\bar{x}_j)|$ and
		\begin{gather}
			\label{ieq:Ebound2}
			|E_\ua^{(k)}(x)|
				\leq C \phi^E_\ua\big((x-\ua),k,b(\ua),b'(\ua),\eps\big),
			\quad
			|E_\oa^{(k)}(x)|
				\leq C \phi^E_\oa\big((\oa-x),k,-b(\oa),b'(\oa),\eps\big).
		\end{gather}		
		The function $\phi^S_{\bar{x}}$ is given as in~\eqref{def:Sphi1} by
		\begin{align}
			\phi^S_{\bar{x}}(x,k,a,b,\eps) & =
			\begin{cases}
				\eps^{\lambda/2} \left(\eps^{1/2}+x\right)^{-\lambda-k},
					& a = 0, \quad b > 0, \\
				\left(\eps^{1/2}+x\right)^{\lambda-k} + \eps \left(\eps^{1/2}+x\right)^{-k-2},
					& a = 0, \quad b < 0, \\
				0,	& \text{otherwise}
			\end{cases} \label{def:Sphi2}
		\intertext{with $0 < \lambda < c(\bar{x})/|b'(\bar{x})|$ whereas the function
		$\phi^E_{\bar{x}}$ is given by} 
			\phi^E_{\bar{x}}(x,k,a,b,\eps) & =
			\begin{cases}
				\eps^{-k} e^{\mu a x/\eps},
					& a < 0, \\
				\eps^{-k/2} e^{-\mu \sqrt{c(\bar{x})} \, x/\sqrt{\eps}},
					& a = b = 0, \\
				0,	& \text{otherwise}
			\end{cases} \label{def:Ephi2}
		\end{align}
		with fixed $\mu \in (0,1)$ which slightly differs from~\eqref{def:Ephi1}.
		Furthermore, the residual part $R$ satisfies
		\begin{gather}
			\label{ieq:R}
			|R^{(k)}(x)| \leq C \eps^{-k}
			\qquad \text{and} \qquad LR = F,
		\end{gather}
		where
		\begin{gather}
			\label{ieq:LR}
			\snorm{F^{(k)}(x)}
				\leq C \phi^E_\ua\big((x-\ua),k,b(\ua),b'(\ua),\eps\big)
					+ C \phi^E_\oa\big((\oa-x),k,-b(\oa),b'(\oa),\eps\big)
				\leq C \eps^{-k}.
		\end{gather}
	\end{theorem}
	The theorem can be proven by combining the arguments that will be provided below.
	Note that the residual term $\eps^q R$ also could be hidden in $S$.
	
	\subsection{The smooth and power-type layer part}
	The part $S$ in the decomposition of Theorem~\ref{th:solDecomp2} can be
	defined as solution of a similar problem on a possibly larger interval.
	
	For $\bar{x} \in \{\ua,\oa\}$ we define $\delta_{\bar{x}}^{\#} \geq 0$ by
	\begin{gather*}
		\delta_{\bar{x}}^{\#} =
		\begin{cases}
			1,	& \bar{x} \in \Bexp, \quad \big(\text{i.e.}, \quad b(\bar{x}) \cdot n(\bar{x}) > 0 \quad \text{or} \quad
					(b(\bar{x}) = b'(\bar{x}) = 0),\big) \\
			0,	& \text{otherwise}.
		\end{cases}
	\end{gather*}
	If necessary we extend the functions $b,c$, and $f$ to the interval
	$I^{\#} := [\ua^{\#},\oa^{\#}] := [\ua - \delta_{\ua}^{\#},\oa+\delta_{\oa}^{\#}]$
	such that $b,c,f \in C^q(I^{\#})$. Additionally, $b,c$ should satisfy
	\begin{itemize}
		\item $b(\bar{x}) \cdot n(\bar{x}) < 0$ for all $\bar{x} \in \{\ua^{\#},\oa^{\#}\} \setminus \{\ua,\oa\}$,
		\item all zeros of $b$ in $I^{\#} \setminus \oI$ are multiple zeros,
		\item $c(x) > 0$ in $I^{\#} \setminus \oI$.
	\end{itemize}
	
	Adapting the a priori estimates~\eqref{ieq:solDerBound}, we see that
	the solution $u^{\#}$ of
	\begin{gather*}
	\begin{aligned}
		- \eps u''(x) + b(x)u'(x) + c(x)u(x) & = f(x) \qquad \text{for } x \in I^{\#}, \\
		u(\ua^{\#})= \nu_{-}, \quad u(\oa^{\#}) & = \nu_+,
	\end{aligned}
	\end{gather*}
	has no exponential boundary layers but possibly interior and power-type boundary
	layers in $\oI$. Hence, for $k = 0, \ldots, q$ it holds
	\begin{gather*}
		\begin{aligned}
			|(u^{\#})^{(k)}(x)|
				& \leq C \bigg( 1 + \phi^S_\ua\big((x-\ua),k,b(\ua),b'(\ua),\eps\big)
					+ \phi^S_\oa\big((\oa-x),k,-b(\oa),b'(\oa),\eps\big) \\
				& \hspace{18em} + \sum_j \left( \eps^{1/2} + |x-\bar{x}_j| \right)^{\lambda_j-k} \bigg),
		\end{aligned}
	\end{gather*}
	with $0 < \lambda_j < c(\bar{x}_j)/|b'(\bar{x}_j)|$ and $\phi^S_{\bar{x}}$ as in~\eqref{def:Sphi2}.

	\subsection{General study of exponential boundary layer corrections}
	\label{app:layerCorrections}
	We will study exponential boundary layer corrections in a general setting now.
	So let $\bar{x} \in \{\ua,\oa\}$. In order to cover the different cases we often
	use ``$\pm$'' and ``$\mp$''. Note that then the upper case is always associated
	with $\bar{x} = \ua$ and the lower case with $\bar{x} = \oa$. Furthermore, generally
	suppose $x \in [\ua,\oa]$.
	
	\subsubsection[\texorpdfstring{Layer of width $\KAO(\eps \log 1/\eps)$}{Layer of width epsilon}]
		{Layer of width $\boldsymbol{\KAO(\eps \log 1/\eps)}$}
	\label{app:sub:epsLayerCor}
	
	First consider the transformation $\xi = \mp \frac{\bar{x} - x}{\eps}$. It follows
	\begin{gather*}
		\xi = \mp \frac{\bar{x} - x}{\eps} \qquad \Longleftrightarrow \qquad
		\mp \xi \eps = \bar{x} - x \qquad \Longleftrightarrow \qquad
		x = \bar{x} \pm \xi \eps.
	\end{gather*}
	This transformation is needed when an exponential layer occurs at the outflow boundary
	(on the left (upper case) or on the right (lower case)) and so we assume that $\sgn(b(\bar{x})) = \mp 1$.
	
	Set $v(x) = \tilde{v}(\xi(x))$ and let $\widetilde{L}$ denote the operator $L$ with
	respect to $\xi$. We have
	\begin{gather*}
		\frac{\mathrm{d} v}{\mathrm{d} x}
			= \frac{\mathrm{d} \tilde{v}}{\mathrm{d} \xi} \frac{\mathrm{d} \xi}{\mathrm{d} x}
			= \pm \eps^{-1} \frac{\mathrm{d} \tilde{v}}{\mathrm{d} \xi},
		\qquad \qquad 
		\frac{\mathrm{d}^2 v}{\mathrm{d} x^2}
			= \eps^{-2} \frac{\mathrm{d}^2 \tilde{v}}{\mathrm{d} \xi^2},
	\end{gather*}
	Therefore, it holds
	\begin{gather*}
		\widetilde{L}\tilde{v}
			= - \eps^{-1} \frac{\mathrm{d}^2 \tilde{v}}{\mathrm{d} \xi^2}
				\pm \eps^{-1} b(\bar{x}\pm \eps \xi) \frac{\mathrm{d} \tilde{v}}{\mathrm{d} \xi}
				+ c(\bar{x}\pm \eps \xi) \tilde{v}.
	\end{gather*}
	We (formally) expand $\tilde{v}$ (the solution of $\widetilde{L}\tilde{v} = 0$), into powers
	of $\eps$, e.g., set
	\begin{gather*}
		\tilde{v}(\xi) = \sum_{i=0}^\infty \eps^i \tilde{v}_i(\xi).
	\end{gather*}
	Additionally, Taylor series expansions (formally) yield
	\begin{gather*}
		b(\bar{x}\pm \eps \xi) = \sum_{j=0}^\infty \frac{(\pm \eps \xi)^j}{j!} b^{(j)}(\bar{x}),
		\qquad \qquad
		c(\bar{x}\pm \eps \xi) = \sum_{j=0}^\infty \frac{(\pm \eps \xi)^j}{j!} c^{(j)}(\bar{x}).
	\end{gather*}
	Sorting the summands of $\widetilde{L}\tilde{v}$ with respect to $\eps$ we obtain
	\begin{align*}
		\widetilde{L}\tilde{v}
		& = - \eps^{-1} \sum_{i=0}^\infty \eps^i \tilde{v}_i''(\xi)
			\pm \eps^{-1} \sum_{j=0}^\infty \frac{(\pm \eps \xi)^j}{j!} b^{(j)}(\bar{x})
				\sum_{i=0}^\infty \eps^i \tilde{v}_i'(\xi)
			+ \sum_{j=0}^\infty \frac{(\pm \eps \xi)^j}{j!} c^{(j)}(\bar{x})
				\sum_{i=0}^\infty \eps^i \tilde{v}_i(\xi) \\
		& = - \sum_{i=-1}^\infty \eps^i \tilde{v}_{i+1}''(\xi)
			\pm \sum_{i=-1}^\infty \sum_{j=0}^\infty \eps^{i+j} \frac{(\pm \xi)^j}{j!} b^{(j)}(\bar{x})
				\tilde{v}_{i+1}'(\xi)
			+ \sum_{i=0}^\infty \sum_{j=0}^\infty \eps^{i+j} \frac{(\pm \xi)^j}{j!} c^{(j)}(\bar{x})
				\tilde{v}_i(\xi) \\
		& = - \sum_{k=-1}^\infty \eps^k \tilde{v}_{k+1}''(\xi)
			\pm \sum_{k=-1}^\infty \eps^k \sum_{j=0}^{k+1} \frac{(\pm \xi)^j}{j!} b^{(j)}(\bar{x})
				\tilde{v}_{k+1-j}'(\xi)
			+ \sum_{k=0}^\infty \eps^k \sum_{j=0}^k \frac{(\pm \xi)^j}{j!} c^{(j)}(\bar{x})
				\tilde{v}_{k-j}(\xi) \\
		& = \sum_{k=-1}^\infty \eps^k
			\Bigg[-\Big(\tilde{v}_{k+1}''(\xi) \mp b(\bar{x}) \tilde{v}_{k+1}'(\xi)\Big) \\
		& \hspace{8em} \pm \sum_{j=0}^k \frac{(\pm \xi)^{j+1}}{(j+1)!} b^{(j+1)}(\bar{x}) \tilde{v}_{k-j}'(\xi)
				+ \sum_{j=0}^k \frac{(\pm \xi)^j}{j!} c^{(j)}(\bar{x}) \tilde{v}_{k-j}(\xi)
			\Bigg] \\
		& = \sum_{k=-1}^\infty \eps^k
			\Bigg[-\Big(\tilde{v}_{k+1}''(\xi) \mp b(\bar{x}) \tilde{v}_{k+1}'(\xi)\Big) \\
		& \hspace{8em} + \sum_{j=0}^k (\pm 1)^j
					\left( \frac{\xi^{j+1}}{(j+1)!} b^{(j+1)}(\bar{x}) \tilde{v}_{k-j}'(\xi)
					+ \frac{\xi^j}{j!} c^{(j)}(\bar{x}) \tilde{v}_{k-j}(\xi) \right)
			\Bigg]\!.
	\end{align*}
	Hence, this results in the following conditions for $\tilde{v}_i$ for $i = 0,1, 2,\ldots$
	\begin{gather*}
		\tilde{v}_{i}''(\xi) \mp b(\bar{x}) \tilde{v}_{i}'(\xi) = \beta_{i-1}(\xi)
			:= \sum_{j=0}^{i-1} (\pm 1)^j
					\left( \frac{\xi^{j+1}}{(j+1)!} b^{(j+1)}(\bar{x}) \tilde{v}_{i-1-j}'(\xi)
					+ \frac{\xi^j}{j!} c^{(j)}(\bar{x}) \tilde{v}_{i-1-j}(\xi) \right)\!.
	\end{gather*}
	In order to guarantee that the correction is local only, we additionally demand that
	\begin{gather*}
		\lim_{\xi \to \infty} \tilde{v}_i(\xi) = 0 \qquad \text{for} \qquad i = 0,1,2,\ldots,
	\end{gather*}
	and in order to enable the ``correct'' boundary conditions
	\begin{gather*}
		\tilde{v}_0(0) = u(\bar{x})-S(\bar{x}), \qquad
		\tilde{v}_i(0) = 0 \qquad \text{for} \qquad i = 1,2,\ldots.
	\end{gather*}
	
	From these conditions we easily get that
	\begin{gather*}
		\tilde{v}_0(\xi) = \left[u(\bar{x})-S(\bar{x})\right] e^{\pm b(\bar{x}) \xi}.
	\end{gather*}
	The method of undetermined coefficients (judicious guessing) inductively implies
	that $\beta_{i-1}$ and $\tilde{v}_i$ are of the form $p(\xi) e^{\pm b(\bar{x}) \xi}$
	for all $i = 0,1,2\ldots$ where $p$ is a suitable polynomial in $\xi$.
	Moreover, following~\cite[p.~94]{O'M91} we even can derive an exact representation
	of $\tilde{v}_i$ by applying variation of constants. So for $i = 1,2,\ldots$ we get
	\begin{gather*}
		\tilde{v}_i(\xi)
			= e^{\pm b(\bar{x})\xi}
			\left(\tilde{v}_i(0)- \int_0^\xi e^{\mp b(\bar{x}) t} \int_t^\infty \beta_{i-1}(s) ds \, dt \right)\!.
	\end{gather*}
	Note that the integral from $(t,\infty)$ always exists since one can recursively
	show that $\beta_{i-1}$ is exponentially decreasing.
	
	Using the well known inequality $1+\xi \leq e^\xi$ and exploiting the known
	structure of $\tilde{v}_i$ we easily get the estimates
	\begin{gather}
		\label{app:ieq:vBoundEps_wrtXi}
		\snorm{\tilde{v}_i^{(k)}(\xi)} \leq C e^{\pm \mu b(\bar{x})\xi},
			\qquad 0 \leq k \leq \tilde{q} \quad \text{for some fixed $\tilde{q} \geq 0$}
	\end{gather}
	with $\mu \in [0,1)$, where $C$ especially depends on $\mu$ and $b(\bar{x})$.
	An inverse transformation yields
	\begin{gather}
		\label{app:ieq:vBoundEps}
		\snorm{v_i^{(k)}(x)} \leq C \eps^{-k} e^{-\mu b(\bar{x})(\bar{x}-x)/\eps},
			\qquad 0 \leq k \leq \tilde{q} \quad \text{for some fixed $\tilde{q} \geq 0$, $\mu \in [0,1)$}.
	\end{gather}
	
	Now, we are interested in $\widetilde{L}\left(\sum_{i=0}^q \eps^i \tilde{v}_i\right)$.
	We assume that $b \in C^{q,1}$ and $c \in C^{q-1,1}$. Then for $q \geq 1$ we have
	\begin{gather*}
		b(\bar{x} \pm \xi \eps)
		= \sum_{j=0}^q \frac{(\pm \xi \eps)^j}{j!} b^{(j)}(\bar{x})
			+ \int_{\bar{x}}^{\bar{x}\pm \xi \eps}
				\frac{(\bar{x}\pm\xi \eps - s)^{q-1}}{(q-1)!} \left[b^{(q)}(s)-b^{(q)}(\bar{x})\right] ds.
	\end{gather*}
	For $q \geq 2$ we have
	\begin{gather*}
		c(\bar{x} \pm \xi \eps)
		= \sum_{j=0}^{q-1} \frac{(\pm \xi \eps)^j}{j!} c^{(j)}(\bar{x})
			+ \int_{\bar{x}}^{\bar{x}\pm \xi \eps}
				\frac{(\bar{x}\pm\xi \eps - s)^{q-2}}{(q-2)!} \left[c^{(q-1)}(s)-c^{(q-1)}(\bar{x})\right] ds
	\end{gather*}
	and for $q \geq 1$
	\begin{gather*}
		c(\bar{x} \pm \xi \eps) = c(\bar{x}) + \left[c(\bar{x} \pm \xi \eps) - c(\bar{x})\right],
	\end{gather*}
	respectively. This yields (for $q \geq 2$)
	\begin{align*}
		& \widetilde{L}\!\left(\sum_{i=0}^q \eps^i \tilde{v}_i\!\right) \\
		& \quad = -\eps^{-1}\sum_{i=0}^q \eps^i \tilde{v}_i''(\xi)
			\pm \eps^{-1} \sum_{j=0}^q \frac{(\pm \xi \eps)^j}{j!} b^{(j)}(\bar{x})
				\sum_{i=0}^q \eps^i \tilde{v}_i'(\xi)
			+ \sum_{j=0}^{q-1} \frac{(\pm \xi \eps)^j}{j!} c^{(j)}(\bar{x})
				\sum_{i=0}^q \eps^i \tilde{v}_i(\xi) \\
		& \hspace{8em} \left.
			\begin{aligned}
				& \pm \eps^{-1} \int_{\bar{x}}^{\bar{x}\pm \xi \eps} \frac{(\bar{x}\pm\xi \eps - s)^{q-1}}{(q-1)!}
					\left[b^{(q)}(s)-b^{(q)}(\bar{x})\right] ds \sum_{i=0}^q \eps^i \tilde{v}_i'(\xi) \\
				& + \int_{\bar{x}}^{\bar{x}\pm \xi \eps} \frac{(\bar{x}\pm\xi \eps - s)^{q-2}}{(q-2)!}
					\left[c^{(q-1)}(s)-c^{(q-1)}(\bar{x})\right] ds \sum_{i=0}^q \eps^i \tilde{v}_i(\xi)
			\end{aligned}\right\} =: \KAI_1 \\
		& \quad = \KAI_1 -\sum_{i=-1}^{q-1} \eps^i \tilde{v}_{i+1}''(\xi)
			\pm \eps^{-1} \sum_{i=0}^q \sum_{j=0}^q \eps^{i+j}
				\frac{(\pm \xi)^j}{j!} b^{(j)}(\bar{x})	\tilde{v}_i'(\xi)
			+ \sum_{i=0}^q \sum_{j=0}^{q-1} \eps^{i+j} \frac{(\pm \xi)^j}{j!} c^{(j)}(\bar{x}) \tilde{v}_i(\xi).
	\end{align*}
	In the ``$b$-sum'' we set $k = i+j$. Note that then $0 \leq k \leq 2q$ and from $0 \leq i = k-j \leq q$
	we additionally have $k-q \leq j \leq k$. Therefore
	\begin{align*}
		\widetilde{L}\!\left(\sum_{i=0}^q \eps^i \tilde{v}_i\!\right)
		& = \KAI_1 -\sum_{k=-1}^{q-1} \eps^k \tilde{v}_{k+1}''(\xi)
			\pm \eps^{-1} \sum_{k=0}^{2q} \eps^{k} \sum_{j=\max\{0,k-q\}}^{\min\{q,k\}}
				\frac{(\pm \xi)^j}{j!} b^{(j)}(\bar{x})	\tilde{v}_{k-j}'(\xi) \\
		& \hspace{12em} + \sum_{k=0}^{2q-1} \eps^{k} \sum_{j=\max\{0,k-q\}}^{\min\{q-1,k\}}
				\frac{(\pm \xi)^j}{j!} c^{(j)}(\bar{x}) \tilde{v}_{k-j}(\xi) \\
		& = \KAI_1 -\sum_{k=-1}^{q-1} \eps^k \tilde{v}_{k+1}''(\xi)
			+ \sum_{k=-1}^{2q-1} \eps^{k} \sum_{j=\max\{-1,k-q\}}^{\min\{q-1,k\}}
				(\pm 1)^j \frac{\xi^{j+1}}{(j+1)!} b^{(j+1)}(\bar{x}) \tilde{v}_{k-j}'(\xi) \\
		& \hspace{12em} + \sum_{k=0}^{2q-1} \eps^{k} \sum_{j=\max\{0,k-q\}}^{\min\{q-1,k\}}
				(\pm 1)^j \frac{\xi^j}{j!} c^{(j)}(\bar{x}) \tilde{v}_{k-j}(\xi).
	\end{align*}
	Now splitting up the sums and noting that the minimums and maximums are known then, we get
	\begin{alignat*}{2}
		\widetilde{L}\!\left(\sum_{i=0}^q \eps^i \tilde{v}_i\!\right)
		& = \KAI_1
			&& + \sum_{k=-1}^{q-1} \eps^k \Bigg[-\Big(\tilde{v}_{k+1}''(\xi) \mp b(\bar{x})\tilde{v}_{k+1}'(\xi)\Big) \\
		&&& \hspace{6em} + \sum_{j=0}^{k}
				(\pm 1)^j \left(\frac{\xi^{j+1}}{(j+1)!} b^{(j+1)}(\bar{x}) \tilde{v}_{k-j}'(\xi)
							+ \frac{\xi^j}{j!} c^{(j)}(\bar{x}) \tilde{v}_{k-j}(\xi) \right) \Bigg]\\
		&&& + \sum_{k=q}^{2q-1} \eps^{k} \Bigg[\sum_{j=k-q}^{q-1}
				(\pm 1)^j \left(\frac{\xi^{j+1}}{(j+1)!} b^{(j+1)}(\bar{x}) \tilde{v}_{k-j}'(\xi)
							+ \frac{\xi^j}{j!} c^{(j)}(\bar{x}) \tilde{v}_{k-j}(\xi) \right) \Bigg].
	\end{alignat*}
	The sum from $k= - 1$ to $q-1$ vanishes due to the construction
	of $\tilde{v}_{k+1}$. Since $b \in C^{q,1}$ and $c \in C^{q-1,1}$,
	we have
	$\left| b^{(q)}(s) - b^{(q)}(\bar{x})\right|
		+ \left| c^{(q-1)}(s) - c^{(q-1)}(\bar{x})\right| \leq C (s-\bar{x})$.
	Hence, using the rule for differentiation under the integral sign we get
	for $0 \leq k \leq q-1$ that
	\begin{gather}
	\label{app:ieq:KAI1_bTerm}
	\begin{aligned}
		& \left|\frac{d^k}{d \xi^k}\int_{\bar{x}}^{\bar{x}\pm \xi \eps} \frac{(\bar{x}\pm\xi \eps - s)^{q-1}}{(q-1)!}
				\left[b^{(q)}(s)-b^{(q)}(\bar{x})\right] ds \right| \\
		& \qquad = \left| (\pm\eps)^k \int_{\bar{x}}^{\bar{x}\pm \xi \eps} \frac{(\bar{x}\pm\xi \eps - s)^{q-1-k}}{(q-1-k)!}
				\left[b^{(q)}(s)-b^{(q)}(\bar{x})\right] ds \right|
			\leq C \eps^{q+1} \xi^{q+1-k}
	\end{aligned}
	\end{gather}
	and analogously for $0 \leq k \leq q-2$ that
	\begin{gather}
		\label{app:ieq:KAI1_cTerm}
		\left|\frac{d^k}{d \xi^k}\int_{\bar{x}}^{\bar{x}\pm \xi \eps} \frac{(\bar{x}\pm\xi \eps - s)^{q-2}}{(q-2)!}
				\left[c^{(q-1)}(s)-c^{(q-1)}(\bar{x})\right] ds \right|
			\leq C \eps^q \xi^{q-k}.
	\end{gather}
	Combining this, Leibniz' rule for the $k$th derivative of a product, \eqref{app:ieq:vBoundEps_wrtXi},
	and the standard inequality $1 + \xi \leq e^{\xi}$, we obtain
	\begin{gather*}
		\left| \frac{d^k}{d \xi^k} \widetilde{L}\!\left(\sum_{i=0}^q \eps^i \tilde{v}_i\!\right) \right|
			\leq C \eps^q e^{\pm \mu b(\bar{x})\xi},
			\qquad 0 \leq k \leq q-2 \quad \text{ with fixed $\mu \in [0,1)$}.
	\end{gather*}
	Noting that $\frac{d^k}{d x^k} {L}\!\left(\sum_{i=0}^q \eps^i v_i\!\right)
	= \frac{d^k}{d \xi^k} \widetilde{L}\!\left(\sum_{i=0}^q \eps^i \tilde{v}_i\!\right) \xi'(x)^{k}$,
	an inverse transformation implies
	\begin{gather}
		\label{app:ieq:LsumBound_Eps}
		\left| \frac{d^k}{d x^k} {L}\!\left(\sum_{i=0}^q \eps^i v_i\!\right) \right|
			\leq C \eps^{q-k} e^{- \mu b(\bar{x})(\bar{x}-x)/\eps},
			\qquad 0 \leq k \leq q-2 \quad \text{with fixed $\mu \in [0,1)$}.
	\end{gather}
	
	\subsubsection[\texorpdfstring{Layer of width $\KAO(\sqrt{\eps} \log 1/\sqrt{\eps})$}{Layer of width root of epsilon}]
		{Layer of width $\boldsymbol{\KAO(\sqrt{\eps} \log 1/\sqrt{\eps})}$}
	\label{app:sub:sqrtEpsLayerCor}
	
	Let us consider the transformation given by $\eta = \mp \frac{\bar{x} - x}{\sqrt{\eps}}$, i.e.,
	\begin{gather*}
		\eta = \mp \frac{\bar{x} - x}{\sqrt{\eps}} \qquad \Longleftrightarrow \qquad
		\mp \eta \sqrt{\eps} = \bar{x} - x \qquad \Longleftrightarrow \qquad
		x = \bar{x} \pm \eta \sqrt{\eps}.
	\end{gather*}
	This transformation is needed when an exponential layer occurs at an multiple boundary
	turning point (on the left (upper) or on the right (lower case)) and so we assume that
	$b(\bar{x}) = b'(\bar{x}) = 0$.
	
	Set $v(x) = \hat{v}(\eta(x))$ and let $\widehat{L}$ denote the operator $L$ with
	respect to $\eta$. Here we have
	\begin{gather*}
		\frac{\mathrm{d} v}{\mathrm{d} x}
			= \frac{\mathrm{d} \hat{v}}{\mathrm{d} \eta} \frac{\mathrm{d} \eta}{\mathrm{d} x}
			= \pm \eps^{-1/2} \frac{\mathrm{d} \hat{v}}{\mathrm{d} \eta},
		\qquad \qquad 
		\frac{\mathrm{d}^2 v}{\mathrm{d} x^2}
			= \eps^{-1} \frac{\mathrm{d}^2 \hat{v}}{\mathrm{d} \eta^2}.
	\end{gather*}
	and it follows
	\begin{gather*}
		\widehat{L}\hat{v}
			= - \frac{\mathrm{d}^2 \hat{v}}{\mathrm{d} \eta^2}
				\pm \eps^{-1/2} b(\bar{x}\pm \sqrt{\eps} \eta) \frac{\mathrm{d} \hat{v}}{\mathrm{d} \eta}
				+ c(\bar{x}\pm \sqrt{\eps} \eta) \hat{v}.
	\end{gather*}
	Expanding the solution $\hat{v}$ of $\widehat{L}\hat{v} = 0$ (formally) into powers
	of $\sqrt{\eps}$ we get
	\begin{gather*}
		\hat{v}(\eta) = \sum_{i=0}^\infty \eps^{i/2} \hat{v}_i(\eta).
	\end{gather*}
	From two Taylor series expansions we (formally) gain
	\begin{gather*}
		b(\bar{x}\pm \sqrt{\eps} \eta) = \sum_{j=0}^\infty \frac{(\pm \sqrt{\eps} \eta)^j}{j!} b^{(j)}(\bar{x}),
		\qquad \qquad
		c(\bar{x}\pm \sqrt{\eps} \eta) = \sum_{j=0}^\infty \frac{(\pm \sqrt{\eps} \eta)^j}{j!} c^{(j)}(\bar{x}).
	\end{gather*}
	Sorting the summands of $\widehat{L}\hat{v}$ with respect to $\sqrt{\eps}$ and using
	$b(\bar{x}) = b'(\bar{x}) = 0$, we obtain
	\begin{align*}
		\widehat{L}\hat{v}
		& = - \sum_{i=0}^\infty \eps^{i/2} \hat{v}_i''(\eta)
			\pm \eps^{-1/2} \sum_{j=0}^\infty \frac{(\pm \sqrt{\eps} \eta)^j}{j!} b^{(j)}(\bar{x})
				\sum_{i=0}^\infty \eps^{i/2} \hat{v}_i'(\eta)
			+ \sum_{j=0}^\infty \frac{(\pm \sqrt{\eps} \eta)^j}{j!} c^{(j)}(\bar{x})
				\sum_{i=0}^\infty \eps^{i/2} \hat{v}_i(\eta) \\
		& = - \sum_{i=0}^\infty \eps^{i/2} \hat{v}_i''(\eta)
			\pm \sum_{i=-1}^\infty \sum_{j=0}^\infty \eps^{(i+j)/2} \frac{(\pm \eta)^j}{j!} b^{(j)}(\bar{x})
				\hat{v}_{i+1}'(\eta)
			+ \sum_{i=0}^\infty \sum_{j=0}^\infty \eps^{(i+j)/2} \frac{(\pm \eta)^j}{j!} c^{(j)}(\bar{x})
				\hat{v}_i(\eta) \\
		& = - \sum_{k=0}^\infty \eps^{k/2} \hat{v}_{k}''(\eta)
			\pm \sum_{k=-1}^\infty \eps^{k/2} \sum_{j=0}^{k+1} \frac{(\pm \eta)^j}{j!} b^{(j)}(\bar{x})
				\hat{v}_{k+1-j}'(\eta)
			+ \sum_{k=0}^\infty \eps^{k/2} \sum_{j=0}^k \frac{(\pm \eta)^j}{j!} c^{(j)}(\bar{x})
				\hat{v}_{k-j}(\eta) \\
		& = \sum_{k=0}^\infty \eps^{k/2}
			\Bigg[-\Big(\hat{v}_{k}''(\eta) - c(\bar{x}) \hat{v}_{k}(\eta)\Big)
				\pm b(\bar{x}) \hat{v}_{k+1}'(\eta)  + \eta b'(\bar{x}) \hat{v}_{k}'(\eta)\\
		& \hspace{7em} \pm \sum_{j=1}^k \frac{(\pm \eta)^{j+1}}{(j+1)!} b^{(j+1)}(\bar{x}) \hat{v}_{k-j}'(\eta)
				+ \sum_{j=1}^k \frac{(\pm \eta)^j}{j!} c^{(j)}(\bar{x}) \hat{v}_{k-j}(\eta)
			\Bigg] \pm \eps^{-1/2} b(\bar{x})\hat{v}_0'(\eta) \\
		& = \sum_{k=0}^\infty \eps^{k/2}
			\Bigg[-\Big(\hat{v}_{k}''(\eta) - c(\bar{x}) \hat{v}_{k}(\eta)\Big) \\
		& \hspace{7em} + \sum_{j=1}^k (\pm 1)^j
					\left( \frac{\eta^{j+1}}{(j+1)!} b^{(j+1)}(\bar{x}) \hat{v}_{k-j}'(\eta)
					+ \frac{\eta^j}{j!} c^{(j)}(\bar{x}) \hat{v}_{k-j}(\eta) \right)
			\Bigg]\!.
	\end{align*}
	From this we get the following conditions for $\hat{v}_i$ with $i = 0,1,2,\ldots$
	\begin{gather*}
		\hat{v}_i''(\eta) - c(\bar{x}) \hat{v}_i(\eta)
		= \gamma_{i-1}(\eta)
		:= \sum_{j=1}^i (\pm 1)^j \left( \frac{\eta^{j+1}}{(j+1)!} b^{(j+1)}(\bar{x}) \hat{v}_{i-j}'(\eta)
					+ \frac{\eta^j}{j!} c^{(j)}(\bar{x}) \hat{v}_{i-j}(\eta) \right)\!.
	\end{gather*}
	
	Since we want the corrections to be local only, we demand that
	\begin{gather*}
		\lim_{\eta \to \infty} \hat{v}_i(\eta) = 0 \qquad \text{for} \qquad i = 0,1,2,\ldots,
	\end{gather*}
	and additionally to enable the ``correct'' boundary conditions
	\begin{gather*}
		\hat{v}_0(0) = u(\bar{x})-S(\bar{x}), \qquad
		\hat{v}_i(0) = 0 \qquad \text{for} \qquad i = 1,2,\ldots.
	\end{gather*}
	Obviously $\gamma_{-1} = 0$. Therefore it is very easy to calculate that
	\begin{gather*}
		\hat{v}_0(\eta) = \left[u(\bar{x})-S(\bar{x})\right] e^{- \sqrt{c(\bar{x})} \eta}.
	\end{gather*}
	Using the method of undetermined coefficients (judicious guessing) we can conclude inductively
	that $\gamma_{i-1}$ and $\hat{v}_i$	are of the form $p(\eta) e^{-\sqrt{c(\bar{x})} \eta}$ for
	all $i = 0,1,2\ldots$ where $p$ is a suitable polynomial in $\eta$. Therefore,
	the inequality $1+\eta \leq e^\eta$ gives
	\begin{gather}
		\label{app:ieq:vBoundSqrtEps_wrtEta}
		\snorm{\hat{v}_i^{(k)}(\eta)} \leq C e^{- \mu \sqrt{c(\bar{x})}\eta},
			\qquad 0 \leq k \leq \hat{q} \quad \text{for some fixed $\hat{q} \geq 0$}
	\end{gather}
	with $\mu \in [0,1)$, where $C$ especially depends on $\mu$ and $\sqrt{c(\bar{x})}$.
	An inverse transformation yields
	\begin{gather}
		\label{app:ieq:vBoundSqrtEps}
		\snorm{v_i^{(k)}(x)} \leq C \eps^{-k/2} e^{\pm \mu \sqrt{c(\bar{x})}(\bar{x}-x)/\sqrt{\eps}},
			\qquad 0 \leq k \leq \hat{q} \quad \text{for some fixed $\hat{q} \geq 0$, $\mu \in [0,1)$}.
	\end{gather}
	
	We now study $\widehat{L}\left(\sum_{i=0}^q \eps^{i/2} \hat{v}_i\right)$.
	Suppose that $b \in C^{q,1}$ and $c \in C^{q-1,1}$. Then for $q \geq 1$ we have
	\begin{gather*}
		b(\bar{x} \pm \eta \sqrt{\eps})
		= \sum_{j=0}^q \frac{(\pm \eta \sqrt{\eps})^j}{j!} b^{(j)}(\bar{x})
			+ \int_{\bar{x}}^{\bar{x}\pm \eta \sqrt{\eps}}
				\frac{(\bar{x}\pm\eta \sqrt{\eps} - s)^{q-1}}{(q-1)!} \left[b^{(q)}(s)-b^{(q)}(\bar{x})\right] ds.
	\end{gather*}
	For $q \geq 2$ we have
	\begin{gather*}
		c(\bar{x} \pm \eta \sqrt{\eps})
		= \sum_{j=0}^{q-1} \frac{(\pm \eta \sqrt{\eps})^j}{j!} c^{(j)}(\bar{x})
			+ \int_{\bar{x}}^{\bar{x}\pm \eta \sqrt{\eps}}
				\frac{(\bar{x}\pm\eta \sqrt{\eps} - s)^{q-2}}{(q-2)!} \left[c^{(q-1)}(s)-c^{(q-1)}(\bar{x})\right] ds
	\end{gather*}
	and for $q \geq 1$
	\begin{gather*}
		c(\bar{x} \pm \eta \sqrt{\eps}) = c(\bar{x}) + \left[c(\bar{x} \pm \eta \sqrt{\eps}) - c(\bar{x})\right],
	\end{gather*}
	respectively. This yields (for $q \geq 2$)
	\begin{align*}
		& \widehat{L}\!\left(\sum_{i=0}^q \eps^{i/2} \hat{v}_i\!\right) \\
		& \quad = -\sum_{i=0}^q \eps^{i/2} \hat{v}_i''(\eta)
			\pm \eps^{-1/2} \sum_{j=0}^q \frac{(\pm \eta \sqrt{\eps})^j}{j!} b^{(j)}(\bar{x})
				\sum_{i=0}^q \eps^{i/2} \hat{v}_i'(\eta)
			+ \sum_{j=0}^{q-1} \frac{(\pm \eta \sqrt{\eps})^j}{j!} c^{(j)}(\bar{x})
				\sum_{i=0}^q \eps^{i/2} \hat{v}_i(\eta) \\
		& \hspace{6em} \left.
			\begin{aligned}
				& \pm \eps^{-1/2} \int_{\bar{x}}^{\bar{x}\pm \eta \sqrt{\eps}} \frac{(\bar{x}\pm\eta \sqrt{\eps} - s)^{q-1}}{(q-1)!}
					\left[b^{(q)}(s)-b^{(q)}(\bar{x})\right] ds \sum_{i=0}^q \eps^{i/2} \hat{v}_i'(\eta) \\
				& + \int_{\bar{x}}^{\bar{x}\pm \eta \sqrt{\eps}} \frac{(\bar{x}\pm\eta \sqrt{\eps} - s)^{q-2}}{(q-2)!}
					\left[c^{(q-1)}(s)-c^{(q-1)}(\bar{x})\right] ds \sum_{i=0}^q \eps^{i/2} \hat{v}_i(\eta)
			\end{aligned}\right\} =: \KAI_2 \\
		& \quad = \KAI_2 -\sum_{i=0}^{q} \eps^{i/2} \hat{v}_{i}''(\eta) 
				\pm \eps^{-1/2} \sum_{i=0}^q \sum_{j=0}^q \eps^{(i+j)/2}
					\frac{(\pm \eta)^j}{j!} b^{(j)}(\bar{x}) \hat{v}_i'(\eta) \\
		& \hspace{20em}
				+ \sum_{i=0}^q \sum_{j=0}^{q-1} \eps^{(i+j)/2} \frac{(\pm \eta)^j}{j!}
					c^{(j)}(\bar{x}) \hat{v}_i(\eta).
	\end{align*}
	In the ``$b$-sum'' we set $k = i+j$. Note that then $0 \leq k \leq 2q$ and from $0 \leq i = k-j \leq q$
	we additionally have $k-q \leq j \leq k$. So
	\begin{align*}
		\widehat{L}\!\left(\sum_{i=0}^q \eps^{i/2} \hat{v}_i\!\right)
		& = \KAI_2 -\sum_{k=0}^{q} \eps^{k/2} \hat{v}_{k}''(\eta)
			\pm \eps^{-1/2} \sum_{k=0}^{2q} \eps^{k/2} \sum_{j=\max\{0,k-q\}}^{\min\{q,k\}}
				\frac{(\pm \eta)^j}{j!} b^{(j)}(\bar{x})	\hat{v}_{k-j}'(\eta) \\
		& \hspace{12em} + \sum_{k=0}^{2q-1} \eps^{k/2} \sum_{j=\max\{0,k-q\}}^{\min\{q-1,k\}}
				\frac{(\pm \eta)^j}{j!} c^{(j)}(\bar{x}) \hat{v}_{k-j}(\eta) \\
		& = \KAI_2 -\sum_{k=0}^{q} \eps^{k/2} \hat{v}_{k}''(\eta)
			+ \sum_{k=-1}^{2q-1} \eps^{k/2} \sum_{j=\max\{-1,k-q\}}^{\min\{q-1,k\}}
				(\pm 1)^j \frac{\eta^{j+1}}{(j+1)!} b^{(j+1)}(\bar{x}) \hat{v}_{k-j}'(\eta) \\
		& \hspace{12em} + \sum_{k=0}^{2q-1} \eps^{k/2} \sum_{j=\max\{0,k-q\}}^{\min\{q-1,k\}}
				(\pm 1)^j \frac{\eta^j}{j!} c^{(j)}(\bar{x}) \hat{v}_{k-j}(\eta).
	\end{align*}
	From splitting of the sums and evaluating the minimums and maximums we obtain
	\begin{alignat*}{2}
		\widehat{L}\!\left(\sum_{i=0}^q \eps^{i/2} \hat{v}_i\!\right)
		& = \KAI_2
			&& - \eps^{q/2} \hat{v}_q''(\eta) \pm \eps^{-1/2} b(\bar{x}) \hat{v}_{0}'(\eta) \\
		&&& + \sum_{k=0}^{q-1} \eps^{k/2} \Bigg[-\Big(\hat{v}_{k}''(\eta) - c(\bar{x})\hat{v}_{k}(\eta)\Big) 
				 \pm b(\bar{x})\hat{v}_{k+1}'(\eta) + \eta b'(\bar{x}) \hat{v}_k'(\eta) \\
		&&& \hspace{6em} + \sum_{j=1}^{k}
				(\pm 1)^j \left(\frac{\eta^{j+1}}{(j+1)!} b^{(j+1)}(\bar{x}) \hat{v}_{k-j}'(\eta)
							+ \frac{\eta^j}{j!} c^{(j)}(\bar{x}) \hat{v}_{k-j}(\eta) \right) \Bigg]\\
		&&& + \sum_{k=q}^{2q-1} \eps^{k/2} \Bigg[\sum_{j=k-q}^{q-1}
				(\pm 1)^j \left(\frac{\eta^{j+1}}{(j+1)!} b^{(j+1)}(\bar{x}) \hat{v}_{k-j}'(\eta)
							+ \frac{\eta^j}{j!} c^{(j)}(\bar{x}) \hat{v}_{k-j}(\eta) \right) \Bigg].
	\end{alignat*}
	
	Because of $b(\bar{x}) = b'(\bar{x}) = 0$ and the definition of $\hat{v}_k$, the
	sum from $k=0$ to $q-1$ vanishes. With arguments similar to those used to
	derive~\eqref{app:ieq:KAI1_bTerm} and~\eqref{app:ieq:KAI1_cTerm} we get
	\begin{align*}
		\left| \frac{d^k}{d\eta^k}\int_{\bar{x}}^{\bar{x}\pm \eta \sqrt{\eps}} \frac{(\bar{x}\pm\eta \sqrt{\eps} - s)^{q-1}}{(q-1)!}
				\left[b^{(q)}(s)-b^{(q)}(\bar{x})\right] ds \right|
			\leq C \eps^{(q+1)/2} \eta^{q+1-k}, \qquad 0\leq k \leq q-1, \\
		\left| \frac{d^k}{d\eta^k}\int_{\bar{x}}^{\bar{x}\pm \eta \sqrt{\eps}} \frac{(\bar{x}\pm\eta \sqrt{\eps} - s)^{q-2}}{(q-2)!}
				\left[c^{(q-1)}(s)-c^{(q-1)}(\bar{x})\right] ds \right|
			\leq C \eps^{q/2} \eta^{q-k}, \qquad 0\leq k \leq q-2.
	\end{align*}
	Thus, from~\eqref{app:ieq:vBoundSqrtEps_wrtEta} and $1 + \eta \leq e^{\eta}$ we gain
	analogously to the last subsection
	\begin{gather*}
		\left| \frac{d^k}{d \eta^k} \widehat{L}\!\left(\sum_{i=0}^q \eps^i \hat{v}_i\!\right) \right|
			\leq C \eps^{q/2} e^{- \mu \sqrt{c(\bar{x})}\eta},
			\qquad 0 \leq k \leq q-2 \quad \text{ with fixed $\mu \in [0,1)$}
	\end{gather*}
	and after an inverse transformation
	\begin{gather}
		\label{app:ieq:LsumBound_SqrtEps}
		\left| \frac{d^k}{d x^k} {L}\!\left(\sum_{i=0}^q \eps^i v_i\!\right) \right|
			\leq C \eps^{(q-k)/2} e^{\pm \mu \sqrt{c(\bar{x})}(\bar{x}-x)/\sqrt{\eps}},
			\qquad 0 \leq k \leq q-2 \quad \text{with fixed $\mu \in [0,1)$}.
	\end{gather}
	
	\subsection{The exponential boundary layer part(s)}
	
	In order to derive the solution decomposition we use a well known technique.
	In general $u^{\#}$ does not satisfy the boundary conditions of problem~\eqref{prob:generalLin}.
	Therefore boundary layer corrections are introduced which then form the exponential
	boundary layer term $E = E_\ua + E_\oa$ of the solution.
	
	These terms are defined as in Appendix~\ref{app:layerCorrections}.
	If there is a layer of width $\KAO(\sqrt{\eps} \log 1/\sqrt{\eps})$, we assume that
	$b \in C^{2q,1}$ and $c \in C^{2q-1,1}$. If there are layers of width $\KAO(\eps \log 1/\eps)$
	only, it suffices to assume that $b \in C^{q,1}$ and $c \in C^{q-1,1}$ to construct
	the corrections.
	
	In order to simplify the notation let us suppose that the
	solution exhibits a layer of width $\KAO(\eps \log 1/\eps)$ at $\oa$ and
	of width $\KAO(\sqrt{\eps} \log 1/\sqrt{\eps})$ at $\ua$. Let $v_i$ denote
	the boundary layer corrections at $\oa$ and $w_j$ at $\ua$, constructed
	according to Appendix~\ref{app:sub:epsLayerCor} and~\ref{app:sub:sqrtEpsLayerCor},
	respectively. Thus for $0 \leq i \leq q$, $0 \leq j \leq 2q$, and $0 \leq k \leq q$
	we have from~\eqref{app:ieq:vBoundEps} and~\eqref{app:ieq:vBoundSqrtEps}
	\begin{gather*}
		\snorm{v_i^{(k)}(x)} \leq C \eps^{-k} e^{-\mu(\oa-x)b(\oa)/\eps},
		\qquad
		\snorm{w_j^{(k)}(x)} \leq C \eps^{-k/2} e^{-\mu\sqrt{c(\ua)}(x-\ua)/\sqrt{\eps}}
	\end{gather*}
	with $\mu \in [0,1)$ fixed.
	
	\subsection{The residual part}
	We now can decompose $u$ as
	\begin{gather*}
		u = \underbrace{u^{\#}}_{=S}
			+ \underbrace{\sum_{i=0}^q \eps^i v_i}_{=E_\oa}
			+ \underbrace{\sum_{j=0}^{2q} \eps^{j/2} w_j}_{=E_\ua}
			+ \eps^q R.
	\end{gather*}
	It remains to bound $R$ and its derivatives. We have
	\begin{gather*}
		R = \eps^{-q} \left(u-u^{\#}- \sum_{i=0}^q \eps^i v_i - \sum_{j=0}^{2q} \eps^{j/2} w_j\right)\!.
	\end{gather*}
	By construction $R$ is (exponentially) small and can be bounded by a
	constant at the boundary. Furthermore, since $L u = f = L u^{\#}$ on $\oI$, we get from
	the known bounds for the derivatives of $L\left( \sum_{i=0}^q \eps^i v_i \right)$ and
	$L\left( \sum_{j=0}^{2q} \eps^{j/2} w_j \right)$, see~\eqref{app:ieq:LsumBound_Eps}
	and~\eqref{app:ieq:LsumBound_SqrtEps}, that
	\begin{gather*}
		L R = F
		\qquad \text{with} \qquad
		\snorm{F^{(k)}} \leq C \eps^{-k} e^{-\mu(\oa-x)b(\oa)/\eps}
			+ C \eps^{-k/2} e^{-\mu\sqrt{c(\ua)}(x-\ua)/\eps} \leq C \eps^{-k},
	\end{gather*}
	where $0 \leq k \leq q-2$ and $\mu \in [0,1)$ fixed. Since we have assumed without loss of generality
	that $c(x) \geq \gamma > 0$, we can apply the maximum principle and get
	\begin{gather*}
		\norm{R}_{\infty} \leq \max\Big\{|R(\ua)|,|R(\oa)|, \gamma^{-1} \max\nolimits_{x \in \oI} |F(x)| \Big\}.
	\end{gather*}
	We follow~\cite[pp.~74,75]{Lis01} to prove estimates for the derivatives of $R$.
	From the mean value theorem there exists a point $x_0 \in (\ua,\oa)$ such
	that
	\begin{gather*}
		\snorm{R'(x_0)} = (\oa-\ua)^{-1} \snorm{R(\oa)-R(\ua)} \leq C.
	\end{gather*}
	Hence, integrating $LR = F$ from $x_0$ to $x$ gives
	\begin{gather*}
		R'(x) = R'(x_0) + \eps^{-1} \left[b(x)R(x)-b(x_0)R(x_0) \right]
			+ \eps^{-1} \int_{x_0}^x \left[ \left(c(s)-b'(s)\right) R(s)-F(s)\right] ds
	\end{gather*}
	and thus $\snorm{R'} \leq C \eps^{-1}$. Differentiating this identity $k-1$ times,
	we get for $2 \leq k \leq q$
	\begin{align*}
		R^{(k)}(x)
			&= \eps^{-1}\sum_{i=0}^{k-1} \tbinom{k-1}{i} b^{(k-1-i)}(x) R^{(i)}(x) \\
			& \qquad + \eps^{-1} \sum_{i=0}^{k-2} \tbinom{k-2}{i}\left( c^{(k-2-i)}(x)-b^{(k-1-i)}(x) \right)R^{(i)}(x)
				- \eps^{-1} F^{(k-2)}(x)
	\end{align*}
	from which we inductively get $\snorm{R^{(k)}} \leq C \eps^{-k}$.
	Summarizing we found out that
	\begin{gather*}
		\snorm{R^{(k)}(x)} \leq C \eps^{-k} \quad \text{for} \quad 0 \leq k \leq q.
	\end{gather*}
	
	\begin{remark}
		If we could/would prove a convenient a priori estimate for problems whose
		right hand sides can depend on $\eps$ as in~\eqref{ieq:LR}, the regularity
		assumptions $b \in C^{q,1}$ and $c \in C^{q-1,1}$ may also could
		suffice when layers of width $\KAO(\sqrt{\eps} \log 1/\sqrt{\eps})$
		are present. 
	\end{remark}
	
\end{appendix}

\bibliographystyle{plain}
	\bibliography{genLinTPPbiblio}
	
\end{document}